\documentclass[a4paper]{amsart}

\usepackage{amssymb,amsmath}
\usepackage{latexsym}
\usepackage{eucal}
\makeatletter
\@addtoreset{equation}{section}\makeatother

\newtheorem{theo}{Theorem}[section]   
\newtheorem{lem}[theo]{Lemma}
\newtheorem{prop}[theo]{Proposition}

\newtheorem{defn}[theo]{Definition}

\def \mrn {{\mathbb R}^n}

\def \mr {{\mathbb R}}

\def \mcs {{\mathcal S}}

\def \mcm {{\mathcal M}}
\def \mca {{\mathcal A}}

\def \mcf {{\mathcal F}}

\def \mcr {{\mathcal R}}

\def \mcc {{\mathbb C}}

\def \mn {{\mathbb N}}

\def \lan {\langle}   
\def \ran {\rangle}   
\def \del {\delta}   

\def \det {\operatorname{det}}
\def \tW {\widetilde{W}}

\def \ox {\bar{X}}

\newcommand{\vol}{\operatorname{vol}}
\newcommand{\dvol}{d\operatorname{vol}}
\newcommand{\ac}{\operatorname{ac}}
\newcommand{\pp}{\operatorname{pp}}

\def \ha{ {\frac{1}{2}}}

\def \oq {\frac{1}{4}}

\def \p {\partial}

\def \rao#1 {\frac{\p}{\p #1} #1}

\newcommand{\mc}{\mathcal}
\newcommand{\rr}{\mathbb{R}}
\newcommand{\nn}{\mathbb{N}}
\newcommand{\cc}{\mathbb{C}}
\newcommand{\hh}{\mathbb{H}}
\newcommand{\zz}{\mathbb{Z}}

\newcommand{\la}{\lambda}
\newcommand{\eps}{\epsilon}

\newcommand{\pl}{\partial}
\newcommand{\x}{\times}

\newcommand{\til}{\widetilde}
\newcommand{\bbar}{\overline}

\newcommand{\cjd}{\rangle}
\newcommand{\cjg}{\langle}

\newcommand{\demi}{\frac{1}{2}}
\newcommand{\ndemi}{\frac{n}{2}}
\newcommand{\tra}{\textrm{Tr}}

\newcommand{\rang}{\textrm{rank}}

\newcommand{\hn}{\textrm{H}_n}
\newcommand{\he}{\textrm{He}}
\newcommand{\rb}{\textrm{rb}}
\newcommand{\lb}{\textrm{lb}}
\newcommand{\frf}{\textrm{ff}}
\newcommand{\trans}{{^t}\!}
\newcommand{\indic}{\operatorname{1\negthinspace l}}

\def\qed{\hfill$\square$}
\input epsf
\begin{document}
\title{Scattering and Inverse Scattering  on ACH Manifolds}
\author{Colin Guillarmou and Ant\^onio S\'a Barreto}
\address{ Laboratoire Jean-Alexandre Dieudonn\'e \\
Universit\'e de Nice, Parc Valrose \\
06108\\
Nice cedex 2 \\
France}
    \email{cguillar@math.unice.fr}
\address{Department of mathematics\\
         Purdue University\\
        150 N. University Street, West-Lafayette\\    
         IN 47907, USA}
     \email{sabarre@math.purdue.edu}

\thanks{  The first author was  partially supported by French ANR grants no. [JC05-52556] and [JC05-46063], and an Australian National University postdoctoral fellowship. Both authors were supported in part by the National Science Foundation Grant DMS 0500788. }

\begin{abstract}
\noindent We study scattering and inverse scattering theories for asymptotically complex hyperbolic manifolds.
We show the existence of the scattering operator 
as a meromorphic family of operators in the Heisenberg calculus on the boundary, which is a contact manifold
with a pseudohermitian structure. Then we define radiation fields as in the real asymptotically hyperbolic case,
and reconstruct the scattering operator from those fields.  As an application we show that the manifold, including its topology and the metric, are determined up to invariants by the scattering matrix at all energies.
\end{abstract}
\maketitle

%\newpage
%\tableofcontents
%\newpage

\section{Introduction}
Scattering theory and inverse problems for {\it real}  asymptotically hyperbolic manifold
have been extensively studied, see for example \cite{G,G1,GRZ,JSB,ME,ME1,SA} and references cited there.  Their complex analogue, the asymptotically complex hyperbolic manifolds, ACH in short,  have not been studied as much.  They were  introduced by Epstein, Melrose and Mendoza \cite{EMM}, 
 and more recently have also been considered by Biquard \cite{BI} and Biquard-Herzlich \cite{BH}.  This class of manifolds
 contains  certain quotients of the  complex hyperbolic space
by discrete groups, as well as smooth pseudo-convex domains in $\cc^{n+1}$ equipped with 
a K\"ahler metric of Bergman type. The purpose of this work is to extend to the complex case several results which are known for real asymptotically hyperbolic manifolds.

Before discussing asymptotically complex hyperbolic manifolds, we recall certain facts about real asymptotically hyperbolic manifolds. An $(n+1)$-dimensional non compact manifold $X$ equipped with a $C^\infty$ Riemannian metric $g$ is called asymptotically hyperbolic if it compactifies into  a $C^\infty$ manifold $\bar{X}$ with boundary 
$\pl\bar{X},$ and if $\rho$ is a defining function of the boundary $\p X,$ $\rho^2g$ is a $C^\infty$ metric which is non-degenerate up to $\p X,$ and moreover  if $|d \rho|_{\rho^2g}=1$ at $\p X.$ It can be shown, see \cite{GRA1}, that $(X,g)$ is asymptotically hyperbolic if and only if there exists 
a diffeomorphism $\psi:[0,\eps)_t\x\pl\bar{X}\to U\subset \bar{X}$ with $\psi(\{0\}\x\pl\bar{X})=\pl\bar{X}$ 
such that
\begin{equation}\label{psig}
\psi^*g=\frac{dt^2+h(t)}{t^2}\end{equation}
where $h(t),$ $t\in [0,\eps)$ is a  $C^\infty$ $1$-parameter family of $C^\infty$ metrics on $\pl\bar{X}.$  
The function $\rho:=\psi_*t$ is a boundary defining function in $\bar{X}$ near $\pl\bar{X},$ 
which can be extended smoothly to $\bar{X}.$ Note that the boundary represents the geometric infinity of 
$X$, as does the sphere $S^n$ for the hyperbolic space $\hh^{n+1}$.   

The spectrum of  $\Delta_g,$ the Laplacian of $(X,g)$ was studied in \cite{MM}; it  consists a finite pure
point spectrum $\sigma_{\pp}(\Delta),$ which is the set of $L^2(X)$ 
eigenvalues, and an absolutely continuous spectrum
$\sigma_{\ac}(\Delta)$ satisfying
\begin{gather*}
\sigma_{\ac}(\Delta)= \left[ n^2/4,\infty\right) \;\
\text{ and } \;\ \sigma_{\pp}(\Delta)\subset 
\left( 0, n^2/4\right).
\end{gather*}

The resolvent 
\begin{gather*}
R(\la)=\left(\Delta_g-\la(n-\la)\right)^{-1},
\end{gather*}
which is a bounded operator in $L^2(X),$ for $\Re(\la)>\frac{n}{2},$  has a finite meromorphic extension to $\cc\setminus \ha(n-\mn_0),$ where $\mn_0$ is the set of non-zero natural numbers,  as a map from $C_0^\infty(X)\to C^\infty(X),$ \cite{MM}.   The poles of
$R(\la)$ are called resonances.  Here, and throughout the paper, we call a family of operators \emph{finite-meromorphic}
if it is meromorphic, i.e. it has a finite Laurent expansion at each point, 
and the rank of the polar part at a pole has finite rank. The \emph{ finite meromorphic } continuation of
$R(\la)$ to the entire complex plane exists
if and only if $h(t)$ has an even Taylor expansion at $t=0$. If $h(t)$ is not even, $R(\la)$ might have essential singularities at the points $\ha(n-\mn_0),$ \cite{G0}.

It has been shown in \cite{GRZ,JSB} that for $\Re(\la)=n/2,$ $\Im(\la)\not=0,$ and any $f\in C^{\infty}(\pl\bar{X})$, 
there exists a unique $u_\la\in C^{\infty}(X)$ satisfying
\begin{gather*}
(\Delta_g-\la(n-\la))u_\la=0 
\end{gather*}
such that near $\pl\bar{X}$
\begin{gather*}
u_\la=\rho^{n-\la}f_-+\rho^{\la}f_++O(\rho^{\ndemi+1}),\quad f_-=f, \quad f_+\in C^{\infty}(\pl\bar{X}).
\end{gather*}

One can use this to define the scattering operator $S(\la),$ for $\Re(\la)=n/2,$ $\Im(\la)\not=0,$ as a generalized Dirichlet-to-Neumann map 
\begin{gather*}
S(\la): C^\infty(\pl\bar{X}) \longrightarrow C^\infty(\pl\bar{X}) \\
f \longmapsto f_+.
\end{gather*}

 Like the Dirichlet-to-Neumann
map, $S(\la)$ is an elliptic pseudo-differential operator, but of order $2\la-n.$ It  
extends meromorphically to $\cc$. The first author recently studied Krein theory
in even dimension manifolds by introducing a generalized determinant of $S(\la)$ and applied it to analyze 
Selberg zeta function for certain quotient of hyperbolic space by discrete groups of isometries, in continuation
of work by Patterson-Perry \cite{PP}. The second author studied inverse problems using $S(\la).$ 
He first proved with Joshi \cite{JSB} that $S(\la)$ for $\la$ fixed determines
the Taylor expansion of $h(t)$ in (\ref{psig}), then more recently he proves in \cite{SA} 
that the map $\la\to S(\la)$ determines the whole manifold up to global isometry.

We give a precise definition of what we call an \emph{asymptotically complex hyperbolic metric}
in Section \ref{sectach}, but we will briefly explain this notion before stating our results. 
We consider a non-compact Riemannian manifold $(X,g)$ that compactifies into $\bar{X}$
smooth with boundary $\pl\bar{X}.$  We assume that the boundary admits a contact form 
$\Theta_0$ and
an almost complex structure $J:\ker\Theta_0\to\ker\Theta_0$ such that $d\Theta_0(.,J.)$ is 
symmetric positive definite on $\ker \Theta_0.$ The associated Reeb vector field $T_0$  is the one which satisfies 
\begin{gather*}
\Theta_0(T_0)=1 \text{  and  } d\Theta_0(T_0,JZ)=0 \text{ for any } Z\in\ker\Theta_0.
\end{gather*}
The parabolic dilation
$M_\rho$ on $T\pl\bar{X}=\ker\Theta_0\oplus \rr T_0$  is defined by
\[M_\rho(V+tT_0)=\rho V+\rho^2tT_0, \quad V\in\ker\Theta_0, t\in\rr.\] 
We say that $(X,g)$ is ACH if there exists a diffeomorphism $\phi:[0,\eps)_\rho\x\pl\bar{X}\to \phi([0,\eps)_\rho\x\pl\bar{X})\subset \bar{X}$ such that $\phi(\{0\}\x\pl\bar{X})
=\pl\bar{X}$ and
\[\phi^*g=\frac{4d\rho^2+d\Theta_0(.,J.)}{\rho^2}+\frac{\Theta_0^2}{\rho^4}+\rho Q_\rho=:\frac{4d\rho^2+h(\rho)}{\rho^2}\]
for some symmetric tensors $Q_\rho$ on $\pl\bar{X}$ satisfying 
$M_\rho^*Q_\rho \in C^{\infty}([0,\eps)_\rho\x\pl\bar{X}, S^{2}(T\pl\bar{X}))$.  
We call such a $\phi$ a \emph{product decomposition} and we say that $g$ is \emph{even at order $2k$} if 
$h^{-1}(\rho)$ has only even
powers in its Taylor expansion at $\rho=0$ at order $2k$, here $h^{-1}(\rho)$ is the 
metric dual to $h(\rho)$ on $T^*\pl\bar{X}$. We will show that the latter is independent of 
$\phi$.
Note that if $\rho$ is any boundary defining function in $\bar{X}$, then $\rho^4g|_{\pl\bar{X}}=e^{4\omega_0}\Theta_0^2$
for some $\omega_0\in C^\infty(\bar{X})$, thus we have more naturally a conformal class of $1$-form 
$[\Theta_0]$ associated to $g$. 
We call the boundary $\pl\bar{X}$ equipped with $(\Theta_0,J)$ a \emph{pseudo-hermitian structure} and its conformal
class $([\Theta_0],J)$ a \emph{conformal pseudo-hermitian structure}.
On such a manifold, one can define the class $\psi_{\Theta_0}^*(\pl\bar{X})$ 
of Heisenberg pseudo-differential operators associated to $\ker \Theta_0$ and its related principal symbol, 
see \cite{EMM2,Po,BG} and Subsection \ref{heisenpseudo} below.
We can define the ``parabolically homogeneous norm'' on $T\pl\bar{X}$: 
\[||V||_{\he}:=\Big(\Theta_0(V)^2+\frac{1}{4}d\Theta_0(V,JV)^2\Big)^{\frac{1}{4}}.\]
and the metric $h_0:=M_\rho^*(\rho^{-2}h(\rho))|_{\rho=0}=d\Theta_0(.,J.)+\Theta^2_0$.\\

As in the real asymptotically hyperbolic manifolds,  the spectrum of the Laplacian $\Delta_g$ of an ACH manifold $(X,g)$ consists of an absolutely continuous part and the  pure point spectrum satisfying
 \begin{gather}
 \sigma_{ac}=\big[(n+1)^2/4,\infty\big) \text{ and } \sigma_{\rm pp}\subset \big(0,(n+1)^2/4\big), \label{specach}
 \end{gather}
 where $\sigma_{pp}$ is {\it a finite set of eigenvalues.}  The resolvent 
\[R(\la):=(\Delta_g-\la(n+1-\la))^{-1} \in\mc{L}(L^2(X)), \quad \Re(\la)>\frac{n+1}{2},\]
is meromorphic, and Epstein, Melrose and Mendoza  \cite{EMM} proved that it has a finite-meromorphic extension 
to $\cc\setminus(\mc{P}_0\cup -\nn_0),$ where $\mc{P}_0:=\frac{n+1/2}{2}-\demi\nn_0,$
 as a family of pseudo-differential operators in a certain calculus. If we assume that $g$ is even at order $2k$,  $k\in\nn_0,$ it will be shown that $\mc{P}_0$ may be replaced by 
 $\mc{P}_k:=\frac{n+1/2-k}{2}-\demi\nn_0.$ 

As in the real case, we use this strong result to show that for any $\la$ with $\Re(\la)=(n+1)/2,$ and
$\Im(\la)\not=0,$ and for any  $f\in C^{\infty}(\pl\bar{X}),$ 
there exists a unique $u_\la\in C^{\infty}(X)$ satisfying
\begin{gather*}
(\Delta_g-\la(n+1-\la))u_\la=0,
\end{gather*}
 such that near $\pl\bar{X}$
\[u_\la=\rho^{n+1-\la}f_-+\rho^{\la}f_++O(\rho^{\frac{n+1}{2}+1}),\quad f_-=f, \quad f_+\in C^{\infty}(\pl\bar{X}).\] 
We then define the scattering operator 
\begin{gather*}
S(\la):C^\infty(\pl\bar{X})\longrightarrow  C^{\infty}(\pl\bar{X}) \\
f \longmapsto f_+.
\end{gather*}

This operator depends on first derivative of $\rho$ at $\pl\bar{X}$, and thus equivalently on the conformal representative 
in $[\Theta_0]$. For another choice $\hat{\rho}=e^{\omega}\rho$ of boundary defining function  
with $\omega\in C^{\infty}(\bar{X})$, we clearly have the conformal covariance
\[\hat{S}(\la)=e^{-2\omega_0}S(\la)e^{2(n+1-\la)\omega_0}, \quad \omega_0:=\omega|_{\pl\bar{X}}.\] 
The structure of the operator $S(\la)$ is established in
\begin{theo}\label{scattering}
Let $(X,g)$ be an ACH manifold even at order $2k$, then the scattering operator $S(\la)$ extends to $\cc\setminus (-\nn_0\cup\mc{P}_k)$ as a meromorphic 
family of conformally covariant operators in the class of Heisenberg pseudodifferential operators $\Psi_{\Theta_0}^{4\la-2(n+1)}(\pl\bar{X})$,  which is unitary on $L^2(\pl\bar{X},{\rm dvol}_{h_0})$ 
when $\Re(\la)=\frac{n+1}{2}$ and $\Im (\la)\not=0.$ The principal symbol of $S(\la)$ is 
\[\sigma_{\rm pr}(S(\la))(\xi)=c_n\frac{2^{2\la+1}\Gamma(\la)^2}{\Gamma(2\la-n-1)}
\mc{F}_{V\to\xi}(||V||_{\rm{He}}^{-4\la})\]
where $c_n\in\cc$ depends only on $n$ and $\mc{F}$ denotes Fourier transform from $T\pl\bar{X}$ to $T^{*}\pl\bar{X}.$
Moreover, $S(\la)$ is finite-meromorphic in $\cc\setminus (-\nn_0\cup \mc{P}_k\cup(n+1-\mc{P}_k))$
and has at most poles of order $1$ at each $\la_k:=\frac{n+1}{2}+\frac{1}{4}k$ with $k\in\nn$,
the residue of which is a Heisenberg differential operator in $\Psi^{k}_{\Theta_0}(\pl\bar{X})$ plus
a finite rank projector appearing if and only if $\la_k(n-1-\la_k)\in\sigma_{pp}(\Delta_g)$. Moreover  
at $\la_{2k}$, we have 
\[\textrm{Res}_{\la_{2k}}S(\la)=\frac{1}{2((k-1)!k!)}\prod_{l=1}^k(-\Delta_b+i(k+1-2l)T_0) \textrm{ mod }\Psi_{\Theta_0}^{2k-1}(\pl\bar{X})\]
where $\Delta_b,T_0$ are the horizontal sublaplacian and the Reeb vector field of $(\pl\bar{X},\Theta_0,J)$.
\end{theo}

We also deduce from \cite{G1} an explicit formula between finite-multiplicity poles of $S(\la)$ (scattering poles) 
and finite multiplicity poles of $R(\la)$ (resonances) in Proposition \ref{multiplicity} and show that essential singularities
for $S(\la)$ and $R(\la)$ can occur at $(n+1/2-\nn_0)/2$ if the metric has no evenness property, see Proposition
\ref{esssing}.

The proof that $S(\la),$ $\Re(\la) =\frac{n+1}{2},$ $\Im(\la)\not=0,$ is a pseudodifferential operator in the Heisenberg calculus is sketched by Melrose in \cite{MEH}. The novelties in this theorem are the computation of the principal symbol of $S(\la),$  its meromorphic continuation, and the analysis of the poles. In the case where the manifold $\bar{X}$ is a strictly  pseudoconvex domain of $\cc^{n+1}$ equipped with an approximate Einstein K\"ahler metric, the relationship between the residues $\textrm{Res}_{\la_{2k}}S(\la)$ and the Gover-Graham operators of \cite{GoGR} is announced in \cite{HPT}.\\

Then we study the scattering theory from a dynamical view point as in the Lax-Phillips theory.  
We define the radiation fields, show that they give unitary translation representations of the wave group which can be used to define the scattering matrix \eqref{defscat} from the wave equation.

The Cauchy problem for the wave equation
\begin{gather}
\begin{gathered}\label{waveeq}
\left(D_t^2-\Delta_g-\frac{(n+1)^2}{4}\right) u(t,m)=0  \text{ in } \mr_+ \times X \\
u(0,m)=f_1(m), \;\ D_t u(0,m)=f_2(z), \;\ f_1, f_2 \in C_0^\infty(X)
\end{gathered}
\end{gather}
has smooth solutions $u\in C^\infty(\mr_+\times X)$, we consider the behavior of
$u$ at infinity along some bicharacteristics and prove
\begin{theo} Let $(\rho,z)\in[0,\eps)\x\pl\bar{X}$  be some coordinates given by product decomposition 
$\phi$ as above. Let $u(t,z)$ be the solution of \eqref{waveeq}, then 
$$v_+(\rho,s,z'):=\rho^{-n-1}u(s- 2\log \rho, \rho, z)\in C^\infty(\mr\times [0,\eps) \times \pl\bar{X}).$$
\end{theo}

We define the \emph{forward radiation field} as the operator 
\begin{gather}
\begin{gathered}
\mcr_+: C_0^\infty(X) \times C_0^\infty(X) \longrightarrow C^\infty(\mr \times \pl\bar{X}), \\
(f_1,f_2) \longmapsto \frac{\p}{\p s} v_+(0,s,z').
\end{gathered}\label{fwdrf}
\end{gather}
Similarly one can how that $v_-(s,\rho,z):=\rho^{-n-1}u(s+\log \rho,\rho,z)$ 
 is smooth on $\rr\x[0,\eps)\x\pl\bar{X}$ and we can define the \emph{backward radiation field} 
by $\mcr_-(f_1,f_2):=\pl_sv_-(0,s,z)$. 

Let $E_{\textrm{ac}}:=\Pi_{\textrm{ac}}(H^1(X)\x L^2(X))$ where $\Pi_\textrm{ac}$ 
is the orthogonal projection from $L^2(X)$ onto the space of absolute continuity of $\Delta_g,$ and let $H^1(X)$
be the first Sobolev space $H^1(X)=\{f\in L^2(X);|df|_g \in  L^2(X)\}$. The space $E_{\rm ac}$ is a Hilbert space
when equipped with the norm 
\[||(\omega_0,\omega_1)||^2_E:=\demi \int_X(|d\omega_0|^2-\frac{(n+1)^2}{4}|\omega_0|^2+|\omega_1|^2)\textrm{ dvol}_g.\]
Then we show 
\begin{theo}
The forward and backward radiation fields $\mc{R}_\pm$ extend to isometric isomorphisms from $E_{\rm ac}$ to
$L^2(\rr\x\pl\bar{X},dr{\rm dvol}_{h_0}).$ Moreover, the map defined by
\begin{equation}\label{defmcs}
\mc{S}:=\mc{R}_+\mc{R}_-^{-1}: L^2(\rr\x\pl\bar{X},dr{\rm dvol}_{h_0})\to L^2(\rr\x\pl\bar{X},dr{\rm dvol}_{h_0})
\end{equation}
is unitary and is a convolution operator in $s$, and conjugating it with Fourier transform in $s$ we have 
\[\mc{F}\mc{S}\mc{F}^{-1}(\la)= -S(\la).\]
\end{theo}
The operator $\mc{S}$ in (\ref{defmcs}) is the dynamical definition of the scattering operator.
Next, using these tools and after proving a localization result for the support of functions $f\in L^2_{\rm ac}(X)$ for which 
$\mc{R}_+(0,f)=0$ in $s\in(-\infty,s_0)$ (a ``support theorem'' in the sense of Helgason, Lax-Phillips), we are able to prove the following result on inverse scattering:
\begin{theo}\label{inverseprob}
Let $(X_1,g_1), (X_2,g_2)$ two $ACH$ manifolds with the same boundary $M:=\pl\bar{X}_1=\pl\bar{X}_2,$ and equipped 
with the same conformal class of contact forms $[\Theta_{0,1}]=[\Theta_{0,2}]$. Let $S_1(\la),S_2(\la)$
be the corresponding scattering operators associated to a conformal representative $\Theta_0\in[\Theta_{0,1}]$.
If $S_1(\la)=S_2(\la)$ on $\{\Re(\la)=\frac{n+1}{2}, \Im \la\not=0\}$, then there exists a diffeomorphism $\Phi:\bar{X}_1\to\bar{X}_2$
such that $\Phi={\rm Id}$ on $M$ and $\Phi^*g_2=g_1$. 
\end{theo}

The method we use is very close to that introduced by the second author \cite{SA} in the asymptotically 
hyperbolic case, which was inspired by the boundary control theory of Belishev \cite{BE}.\\ 

The paper is organized as follows: In Section 1, we consider the model case of the complex hyperbolic space $\hh^{n+1}_\cc$, then
we discuss the geometry of ACH manifolds near infinity in Section 2. We define
the $\Theta$-calculus on $X$ and the Heisenberg calculus on $\pl\bar{X}$ (these are 
the ``natural'' class of pseudo-differential operators associated to the geometric structure)
in Section 4 and we analyze the Poisson and the scattering operators in Section 6.
The next sections consist in defining raidiation fields (Section 7), prove their relation with scattering 
operator (Section 8), the support Theorem (Section 10) and the inverse problem (Section 11).
We conclude with a technical appendix.\\

\textbf{Acknowledgement} We thank D. Geller, R. Graham, P. Greiner, M. Olbrich and R. Ponge 
for helpful discussions. 

\section{The model case of $\hh^{n+1}_\cc$}\label{modelcase}
\subsection{$\hh^{n+1}_\cc$ and the Heisenberg group $\hn$}
The hyperbolic complex space of complex dimension $n+1$ is denoted by 
$\hh^{n+1}_\cc$, it is the unit ball $B^{n+1}=\{z\in\cc^{n+1}; |z|<1\}$
equipped with the K\"ahler metric 
$g_0:=-4\pl\bar{\pl}\log(\rho)$ where $\rho:=1-|z|^2$.
Note that $\rho$ is a boundary defining function of the 
closed complex ball. The holomorphic curvature is $-1$ 
and this metric is called the \emph{Bergman metric}.
Another model of $\hh^{n+1}_\cc$ is given by 
$\{z\in\cc^{n+1}; Q(z,z)>0\}$ where $Q$ is the quadratic form
\begin{equation}\label{q}
Q(z,z')=-\frac{i}{2}(z_1-\bar{z}_1')-\demi\sum_{j>1}z_j\bar{z}_j',
\end{equation}
and the boundary (the sphere $S^{2n+1}$) is a compactification of the Heisenberg group
\[\hn:=\{z\in\cc^{n+1};Q(z,z)=0\}=\{(\Re(z_1),\demi|\omega|^2,\omega); (z_1,\omega)\in\cc^{n+1}\}
\simeq \rr\x\cc^n\simeq \rr^{2n+1},\] 
thus $\hh^{n+1}_\cc\simeq (0,\infty)\x \hn$. 
The variable $u:=\Re(z_1)$ is the one lying in $\rr$ and we have a contact 
form on $\hn$ given by 
\[\Phi:=du+y.dx-x.dy,\] 
where $\omega=x+iy\in\cc^{n}=\rr^{n}+i\rr^{n}$. The functions  
\[\rho_0:=Q(z,z)^{\demi}, \quad u=\Re(z_1),\quad \omega\in\cc^n\]
give coordinates on $(0,\infty)\x \hn\simeq\hh^{n+1}_\cc$ 
and the Bergman metric with holomorphic curvature $-1$ is given in this model by 
\begin{equation}\label{metricberg}
g_0=\frac{4 d\rho_0^2+2|d\omega|^2}{\rho_0^2}+\frac{\Phi^2}{\rho_0^4}.
\end{equation}
The Heisenberg group $\hn$ is a Lie group with the law
\[(u,\omega)._{\hn} (u',\omega'):=(u+u'-\Im(\omega.\bar{\omega}'),\omega+\omega'),\] 
the origin is $0$ and the inverse $(u,\omega)^{-1}=(-u,-\omega)$.
A basis for the Lie algebra $\textbf{h}$ of $\hn$ is given 
by the left invariant vector fields 
\begin{equation}\label{basehn}
X_j=\frac{1}{\sqrt{2}}(\pl_{x_j}-y_j\pl_u), \quad Y_j=\frac{1}{\sqrt{2}}(\pl_{y_j}+x_j\pl_u), 
\quad T=\pl_u.
\end{equation} 
The map $(u,\omega)\in \rr\x\cc^{n}\mapsto uT+\sum_j\Re(\omega_j)X_j+\Im(\omega_j)Y_j$ 
identifies $\hn$ with $\textbf{h,}$ and the group law becomes
\[W._{\hn} W'=(\Phi(W+W')-d\Phi(W,W'))T+\pi_{\ker\Phi}(W+W'),\]
where $\pi_{\ker\Phi}$ is the projection on $\ker\Phi$ parallel to $T$.

The complex hyperbolic space has a Lie group structure, this is actually a semi-direct product 
of the multiplicative group $((0,\infty),\x)$ with $(\hn,._{\hn})$. We introduce the parabolic
dilation  $M_\delta(\rho_0,u,\omega):=(\delta\rho_0, \delta^2u,\delta\omega)$
on $(0,\infty)\x \hn$ (here $\delta>0$), 
then the group law on $\hh^{n+1}_\cc\simeq (0,\infty)\x \hn$ is
\begin{equation}\label{grouplawhyp}
(\rho_0, W)._{\hh^{n+1}_\cc}(\rho_0',W'):=(\rho_0\rho_0', W._{\hn} M_{\rho_0}(W')).\end{equation}
and we have for this law $(\rho_0,W)^{-1}=(\rho_0^{-1},-M_{\rho_0^{-1}}W)$.  
It is easy to check that the corresponding Lie algebra has a basis 
\begin{equation}\label{basehyp}
\rho_0\pl_{\rho_0}, \rho_0^2\pl_u,\rho_0 X_1,\dots,\rho_0 X_n,\rho_0 Y_1,\dots,\rho_0 Y_n, \end{equation} 
 which is orthonormal with respect to the metric $g_0.$  This algebra will be denoted by 
$^\Phi T_0\hh^{n+1}_\cc,$ to agree with the notation used in the next sections. Observe also that these vectors and the metric
$g_0$ are homogeneous of degree $0$ under the parabolic dilation $M_\delta$.
 
\subsection{The Resolvent kernel for $\hh^{n+1}_\cc$}\label{resolventonhn} 
The spectrum of the Bergman Laplacian $\Delta_{g_0}$ of $\hh^{n+1}_\cc$ is absolutely 
continuous and equal to $\sigma(\Delta_{g_0})=[\frac{(n+1)^2}{4},\infty)$, this leads to study 
the modified resolvent
\[R(s):=(\Delta_{g_0}-s(n+1-s))^{-1}\]
which is bounded on $L^2(\hh^{n+1}_\cc, \textrm{dvol}_{g_0}),$ provided $\Re(s)>\frac{n+1}{2}$.
The Schwartz kernel of $R(s)$ has been computed by Epstein-Melrose-Mendoza \cite{EMM} 
and admits a meromorphic continuation to $\cc$, with poles at $-\nn_0$
of finite multiplicity (contrary to what is written in \cite{EMM}).
By symmetry arguments, this kernel $R(s;z;z')$ is expressed as a 
function of the Bergman distance of $d_{g_0}(z;z')$. We have 
\[\cosh\left(\frac{d_{g_0}(z;z')}{2}\right)=\frac{|Q(z,z')|}{(Q(z,z)Q(z',z'))^\demi}\] 
where $Q$ is defined in (\ref{q}). Using a polar decomposition around the diagonal, the 
kernel $R(s;z,z')$ is obtained as a solution of an hypergeometric ODE, exactly like 
in the real case, and is given by
\begin{gather}
\begin{gathered}
R(s;z,z')=c_n\frac{\Gamma(s)^2}{\Gamma(2s-n)}r(z;z')^s {_2F}_1(s,s,2s-n;r(z;z')), \text{ where } \\
r(z;z'):=\left(\cosh\left(\frac{d_{g_0}(z;z')}{2}\right)\right)^{-2}=\frac{4\rho_0^2{\rho_0'}^2}
{(u-u'+\Im(\omega.\bar{\omega}'))^2+(\rho_0^2+{\rho_0'}^2
+\demi|\omega-\omega'|^2)^2},
\end{gathered}\label{resolvent}
\end{gather}
with $c_n$ constant depending on $n$ and $\textrm{ }_2F_1$ is a hypergeometric function (see \cite{AS}),
we also used the formula
\[Q(z,z')=\frac{-i}{2}(u-u'+\Im(\omega.\bar{\omega}'))+\demi\Big(\rho^2_0(z)+{\rho_0'}^2(z)
+\demi|\omega-\omega'|^2\Big).\]
A change of variables shows that if an operator $K$ has 
a distributional Schwartz kernel which is of the form $k(r(z,z')),$ in other words, it
depends only on $d_{g_0}(z,z'),$  then $K$ is a convolution operator with respect to the group law on $\hh_\cc^{n+1}$:
\begin{equation}\label{convolution}
Kf(\rho_0,u,\omega)=\int k\left(\frac{4\mu^2}{t^2+(1+\mu^2+\demi|z|^2)^2}\right)
f\Big((\rho_0,u,\omega)._{\hh^{n+1}_\cc}(\mu,t,z)^{-1}\Big)\frac{2^{n+1}d\mu dt dz}{\mu},
\end{equation}
where $\mu^{-1}d\mu dtdz$ is a right invariant mesure.  The resolvent kernel \eqref{resolvent} is of this form
($s$ is a parameter), so the action of the operator $R(s)$ on a function is given by \eqref{convolution}.\\

\textsl{Remark}: We see that the poles at $-m\in-\nn_0$ have residue 
\[P_m=\sum_{k=0}^m a_{m,k}r^{k-m}=(Q(z,z)Q(z',z'))^{-m}\sum_{k=0}^m a_{m,k}
|Q(z,z')|^{2m-2k}(Q(z,z)Q(z',z'))^k\]
for some $a_{m,k}\in\cc$. But clearly $P_m$ has finite rank since it is a polynomial 
times $Q(z,z)^{-m}Q(z',z')^{-m}$. So the poles are of finite multiplicity.\\

\section{Asymptotically complex hyperbolic manifolds}\label{sectach}

\subsection{$\Theta$ metrics}
We start by describing the $\Theta$ structures of Epstein-Melrose-Mendoza \cite{EMM}, 
which generalize Bergman type metrics on pseudoconvex domains, as well as 
quotient $\Gamma\backslash\hh^{n+1}_\cc$ of $\hh^{n+1}_\cc$ by convex co-compact groups
of isometries.

Let $\bar{X}=X\cup \p X$ be a smooth $2n+2$-dimensional compact manifold with boundary 
$\pl\bar{X}$ and let $\Theta\in 
C^\infty(\pl\bar{X},T^*\bar{X})$ be a smooth $1$-form on $\pl\bar{X}$ such that 
if $i:\pl\bar{X}\to \bar{X}$ is 
the inclusion, then $\Theta_0:=i^*\Theta$ does not vanish on $\pl\bar{X}.$ 
According to the terminology of \cite{Po}, the boundary $\pl\bar{X}$ has the sturcture of a Heisenberg 
manifold equipped with 
the subbundle $\ker{\Theta_0}$. 

We first recall a few definitions introduced in \cite{EMM}.
If $\rho$ is a boundary defining function of $\pl\bar{X}$, 
we define the Lie subalgebra $\mc{V}_\Theta$ of $C^\infty(\bar{X},T\bar{X})$ 
by the condition 
\[V\in \mc{V}_\Theta \iff V\in\rho C^\infty(\bar{X}, T\bar{X}), 
\til{\Theta}(V) \in \rho^2C^\infty(\bar{X}),\] 
where $\til{\Theta}\in C^\infty(\bar{X},T^*\bar{X})$ is any smooth extension 
of $\Theta$. It is shown in \cite{EMM} that $\mc{V}_\Theta$ only depends 
on the conformal class of $\Theta$.
Let $T,N,Y_1,\dots,Y_{2n}$ be a smooth local frame in $\bar{X}$ near a point 
$p\in\pl\bar{X}$ such that 
\[\textrm{Span}(N,Y_1,\dots,Y_{2n})\subset\ker\til{\Theta}, \quad 
\textrm{Span}(T,Y_1,\dots,Y_{2n})\subset T\pl\bar{X}, \quad d\rho(N)=\til{\Theta}(T)=1.\]
Then any $V\in\mc{V}_\Theta$ can be written near $p$ as 
\begin{equation}\label{formebase}
V=a\rho N+b\rho^2T+\sum_{i=1}^{2n}c_i\rho Y_i, \quad a,b,c_i\in C^\infty(\bar{X})
\end{equation}
and 
\begin{equation}\label{basis}
\rho N,\rho^2T,\rho Y_1,\dots, \rho Y_{2n}
\end{equation} 
form a basis of $\mc{V}_\Theta$ 
over $C^\infty(\bar{X})$ near $p$.
The Lie algebra $\mc{V}_\Theta$ is the set of smooth sections of a vector bundle 
over $T\bar{X}$, we denote by ${^\Theta T}\bar{X}$ this bundle. Let
$F_p$ be the set of vector fields vanishing at $p$ 
if $p\in X$ or the set of vector fields of the form (\ref{formebase}) satisfying
$a(p)=b(p)=c_i(p)=0$ if $p\in\pl\bar{X}$. The fibre 
$^\Theta T_p\bar{X}$ at $p\in\bar{X}$ can be defined by 
$^\Theta T_p\bar{X}:=\mc{V}_\Theta/F_p$. If $p\in\pl\bar{X}$, $^\Theta T_p\bar{X}$
is a Lie algebra, and any vector $v\in{^\Theta T}_p\bar{X}$ can be represented as
\[v=a\rho N+b\rho^2 T+\sum_{i=1}^{2n}c_i\rho Y_i , \quad a,b,c_i\in\rr.\]
The dual bundle $^\Theta T^*\bar{X}$ of $^\Theta T\bar{X}$ has for local 
basis near $p\in\pl\bar{X}$ the dual basis to (\ref{basis}) 
\[ \frac{d\rho}{\rho}, \frac{\til{\Theta}}{\rho^2}, \frac{\alpha_1}{\rho}, \dots, 
\frac{\alpha_2}{\rho}.\] 

A $\Theta$-metric is a smooth positive symmetric $2$-tensor on $^\Theta T^*\bar{X}$ 
\[g\in C^{\infty}(\bar{X},S_+^2(^\Theta T^*\bar{X})).\]
We are interested in the special cases of $\Theta$-metrics for which Epstein, Melrose and Mendoza 
\cite{EMM} proved the meromorphic extension of the resolvent.
We begin by the first assumption, which allows to find particular boundary 
defining functions.

\subsection{Model boundary defining functions}\label{modelbdfs}
Let $g$ be a $\Theta$-metric, it thus restricts to a standard
metric in the interior $X$.
If $\rho$ is a boundary defining function, 
we can define the vector field
$X_\rho\in \mc{V}_\Theta$ as the dual of $d\rho/\rho$ via the metric $g$, i.e.
$g(X_\rho,v)=\rho^{-1}d\rho(v)$ for any $v\in {^\Theta T}\bar{X}$, this is a 
smooth non-vanishing section of $^\Theta T\bar{X}$.
It is clear that in $X$, we have 
\[\nabla^{\rho^2g}\rho =\frac{\nabla^g\rho}{\rho^2}=\frac{X_\rho}{\rho}\]
which extends to a non-vanishing vector in $C^{\infty}(\bar{X},T\bar{X})$
transverse to $\pl\bar{X}$ since $d\rho(X_\rho/\rho)=|d\rho/\rho|_g^2\not=0$ on $\pl\bar{X}$ 
(here $\nabla$ means gradient). 
We first assume that 
\[(\textrm{H1})\quad\quad |X_\rho|_g=\demi \textrm{ on }\pl\bar{X}\]
and it is easy to check that this condition does not depend on $\rho$.
The restriction $\rho^4g|_{T\pl\bar{X}}$ is conformal to the tensor
$\Theta_0^2$, which leads to the definition of the conformal class $[\Theta_0]$.

\begin{lem}\label{bdf}
Let $e^{2\omega_0}\Theta_0\in[\Theta_0]$ with $\omega_0\in C^{\infty}(\pl\bar{X})$, 
then there exists a unique, up to $C_0^\infty(X)$, 
boundary defining function $\rho$ of $\pl\bar{X}$ such that 
$|X_\rho|_g=1/2$ in a neighbourhood of $\pl\bar{X}$ and 
$\rho^4g|_{T\pl\bar{X}}=e^{4\omega_0}\Theta_0^2$.
\end{lem}
\textsl{Proof}:  If $x$ is a boundary defining function we search for a function 
$\omega\in C^\infty(\bar{X})$ such that $\rho:=e^{\omega}x$ satisfies
$|X_\rho|_g^2=|d\rho/\rho|^2_g=1/4$ near $\pl\bar{X}$, this can be rewritten under the form
\begin{gather}
2\frac{X_x(\omega)}{x}+\frac{|d\omega|^2_{g}}{x}=\frac{1/4-|X_x|^2_g}{x}.\label{defvf}
\end{gather}
This a first order non-linear PDE with smooth coefficients since
$|d\omega|^2_{g}=O(x^2)$, $|X_x|_g=1/2$ at $x=0$ and since $X_x/x$ is transverse 
to $\pl\bar{X}$ , it is easy to check that this equation is non-characteristic.
By prescribing the value $\omega|_{x=0}=\omega_0$, we obtain a unique solution
in a neighbourhood of $\pl\bar{X}$.
\qed\\

Such a boundary defining function will be called \emph{a model boundary defining function}.
Let $\phi_t$ be the flow of the vector field $4X_\rho/\rho$, we 
consider the diffeomorphism
\[\phi:\left\{\begin{array}{ccc}
[0,\eps)\x \pl\bar{X} & \to & \phi([0,\eps)\x \pl\bar{X})\subset \bar{X}\\
(t, y) & \to & \phi_t(y)
\end{array}\right.\]
Then $\rho(\phi_t(y))=t$ and for any $Z\in T\pl\bar{X}$
\[\phi^*g(\pl_t,\pl_t)=\frac{4}{t^2}, \quad 
\phi^*g(\pl_t,Z)=\frac{dt(Z)}{t^2}=0.\] 
We will write $t=\rho$ and $X_\rho=\rho\pl_\rho$ for what follows and we call
this diffeomorphism a \emph{product decomposition} near $\pl\bar{X}$.
Note also that $\Theta(\pl_\rho)=0$ since $\rho\pl_\rho\in\mc{V}_\Theta$.
With this product decomposition, the metric $g$ has the form 
\begin{equation}\label{modform1}
g=\frac{4d\rho^2+h(\rho)}{\rho^2}
\end{equation} 
in $(0,\eps)_\rho\x\pl\bar{X}$ with $h(\rho)$ a family of metrics on $T\pl\bar{X}$ for $\rho\not=0$
and such that $h(\rho)/\rho^2\in C^\infty(\pl\bar{X},S^2(^\Theta T\bar{X}))$ depending 
smoothly on $\rho\in[0,\eps)$.

We will say that the metric is \emph{even} is $h(\rho)^{-1}$, as a metric on $T^*\pl\bar{X}$  
has an even Taylor expansion at $\rho=0$ in the product decomposition. It is straightforward to
see that this condition is invariant with respect to the choice of model boundary defining function 
$\rho$ (i.e. of product decomposition), for instance from the proof of Lemma 2.1 in \cite{G0} where the PDE is replaced in our case by (\ref{defvf}).
Indeed, if $x$ is a model boundary defining function and $\rho=e^\omega x$
another one, $\omega$ has to satisfy (\ref{defvf}), that is
\[2\pl_x\omega+x\Big( (\pl_x\omega)^2+|d\omega|^2_{h(x)} \Big)=0\]
and the evenness of the Taylor expansion of $|d\omega|_{h(x)}^2$ at $x=0$ was all that 
we needed in \cite{G0}. Note that \emph{evenness at order $2k$} can also be defined
invariantly by requiring $\pl_\rho^{2j+1}h^{-1}(0)=0$ for all $j<k$ (see again \cite{G0} for 
similar definition in the real case).

\subsection{Additional assumptions}
Following \cite{EMM}, we define for $p\in\pl\bar{X}$ 
the one-dimensional subspace of $^\Theta T_p\bar{X}$
\[K_{2,p}:=\{V\in \rho^2C^\infty(\bar{X},T\bar{X})\}/F_p\]
and the $2n$ dimensional subspace 
\[K_{1,p}:=\{V\in \mc{V}_\Theta; V=\rho W, W \textrm{ tangent to }\pl\bar{X}\}/F_p.\] 
The subspace $K_{1,p}$ is a two-steps nilpotent Lie algebra which is the fibre over $p$ 
of the tangent Lie bundle defined in \cite{Po}. We denote by $K_1,K_2$ the bundles
over $\pl\bar{X}$ whose fibre at $p$ are $K_{1,p},K_{2,p}$.
Near $p\in\pl\bar{X}$, let $(Y_1,\dots,Y_{2n})$ be a local basis 
of $\ker \Theta_0\subset T\pl\bar{X}$ and $T\in T\pl\bar{X}$ 
such that $\Theta_0(T)=1$,
this give a local basis of $T\pl\bar{X}$. 
A basis of $K_{2,p}$ is given by the class of $\rho^2T$ mod $F_p$, whereas
$(\rho Y_1,\dots,\rho Y_{2n},\rho^2 T)$ mod $F_p$ gives a basis 
of $K_{1,p}$. This easily shows that $K_2$ is included in the centre of $K_1$. 

Let us denote $\til{K}_1=\ker (\rho^{-1}d\rho)$ the subbundle of $^\Theta T\bar{X}$,
it is isomorphic to $T\pl\bar{X}$ over $\rho\not=0$ and equal to $K_1$
over $\rho=0$. 
Thus the choice of a function $\rho$ (or product decomposition of $\pl\bar{X}$) 
induces orthogonal decompositions for $g$ (outside $\{\rho=0\}$ for the first one)
\[T\bar{X}\simeq \rr\pl_\rho\oplus T\pl\bar{X},\quad 
^\Theta T\bar{X}\simeq \rr\rho\pl_\rho\oplus \til{K}_1.\]
Using this decomposition, we extend $\Theta_0$
on $\rr\pl_\rho\oplus T\pl\bar{X}$ to be constant with respect to 
$\rho,$ and such that $\Theta_0(\pl_\rho)=0$, in particular $\Theta_0$ is extended 
by $\phi^*\Theta$ at $\{\rho=0\}$. 
Then $\rho^{-2}\Theta_0$ is a smooth section of $^\Theta T^*\bar{X}$ 
and $\ker\Theta_0,\ker(\rho^{-2}\Theta_0)$ are respective subbundle of 
$T\pl\bar{X}, \til{K}_1$.  
We have an isomorphism of vector bundles
\[ \psi : \left\{\begin{array}{ccc}
(T\pl\bar{X}/\ker\Theta_0)\oplus \ker\Theta_0 &\to &K_1\\
(p; T_p\oplus Y_p)& \to & (p;\rho Y+\rho^2T\textrm{ mod }F_p)  
\end{array}\right.,\] 
where $Y,T$ are smooth local sections of $T([0,\eps)_\rho\x\pl\bar{X})$, 
constant with respect to $\rho$, such that $Y\in\ker\Theta_0$, 
$Y(p)=Y_p$, and $T(p)=T_p$ mod $\ker\Theta_0$.
Via $\psi_*$, the form $\Theta_0$ on $(T\pl\bar{X}/\ker\Theta_0)\oplus \ker\Theta_0$ 
is mapped onto the form $\rho^{-2}\Theta_0$ on $K_1$.
The subbundle $(T\pl\bar{X}/\ker\Theta_0)$ is mapped onto $K_2$ by $\psi$ and
$\ker\Theta_0$ onto the bundle $\ker(\rho^{-2}\Theta_0)$.
Similarly the $2$-form $d\Theta_0|_{\ker\Theta_0}$ is mapped onto 
$(\rho^{-2}d\Theta_0)|_{\ker(\rho^{-2}\Theta_0)}$. A local choice of vector $T$ 
transversal to $\ker\Theta_0$ in $T\pl\bar{X}$ 
in a neighbourhood $U_p$ of $p\in\pl\bar{X}$  
fixes a vector $\rho^2 T$ transversal to $\ker(\rho^{-2}\Theta_0)$, 
thus a representative vector $\rho^2T|_{\pl\bar{X}}$ of $K_2$,
and a local basis $\rho^{-1}\alpha_1,\dots,\rho^{-1}\alpha_{2n}$ 
for the annihilator $(\ker(\rho^{-2}\Theta_0))^*$ of $\rr\rho^2T$ in $\til{K}_1^*$ 
(i.e. the dual of $\ker(\rho^{-2}\Theta_0)$) can be chosen.

In view of this discussion, we have that $\rho^{-2}h(\rho)\in C^\infty
([0,\eps)\x\pl\bar{X},S^2_+(\til{K}_1^*))$  
and we can write near a point $p\in\pl\bar{X}$ 
\begin{equation}\label{hrho}
\frac{h(\rho)}{\rho^2}=a\frac{\Theta_0^2}{\rho^4}+
\sum_{i,j=1}^{2n}c_{ij}\frac{\alpha_i\otimes \alpha_j}{\rho^2}+
\sum_{i=1}^{2n}b_i\frac{\alpha_i\otimes\Theta_0}{\rho^3}
\end{equation}
for some functions $a,b_i,c_{ij}\in C^\infty(\bar{X})$. 
Note also that $a|_{\rho=0}$ is globally defined and can be taken to be $1$ 
by changing the conformal representative of $[\Theta_0]$.
 
Let us denote by $g_p$ the metric on $^\Theta T_p\bar{X}$, in terms of 
(\ref{modform1}) and (\ref{hrho}), this is
\begin{equation}\label{modform2}
g_p=4\frac{d\rho^2}{\rho^2}+\frac{\Theta_0^2}{\rho^4}+
\sum_{i,j=1}^{2n}c_{ij}(p)\frac{\alpha_i\otimes\alpha_j}{\rho^2}+
\sum_{i=1}^{2n}b_i(p)\frac{\alpha_i\otimes\Theta_0}{\rho^3}.\end{equation}
The assumptions of \cite{EMM} on the metric correspond to 
the minimal assumptions for which $g_p$ is isometric to the complex 
hyperbolic metric. 
The second assumption made in \cite{EMM} is that 
\[(\textrm{H2})\quad\quad \Theta_0 \textrm{ is a contact form on }\pl\bar{X}\]
which means that $d\Theta_0$ is non-degenerate on $\ker\Theta_0$.
The next hypothesis is that for the orthogonal 
decomposition $K_{1,p}=K_{2,p}\oplus L_p$ for $g_p$, 
the map
\[\mu:\left\{\begin{array}{ccc}
L_p & \to & ^\Theta T_p\bar{X}\\
Z & \to & [\rho\pl_\rho,Z]
\end{array}\right.\]
is the identity. 
Since $L_p$ is spanned by some $Z_i=\sum_{j}l_{ij}\rho Y_j+k_i\rho^2T$ mod $F_p$  
the assumption clearly reduces to $k_i=0$ since $\rho^2T$ mod $F_p$ commutes 
with any elements of $K_{1,p}$, therefore $L_p=(\ker(\Theta_0/\rho^2))_p$. 
Then by orthogonality of the decomposition, 
this means that $b_i(p)=0$, i.e. 
\[(\textrm{H3})\quad \quad \ker\Big(\frac{\Theta_0}{\rho^2}\Big)\perp_g K_2.\]  
The last assumption of \cite{EMM} is 
\[(\textrm{H4})\quad \quad \exists \; J\in \textrm{End}
\Big(\ker\Big(\frac{\Theta_0}{\rho^2}\Big)\Big), 
J^2=-\textrm{Id} \textrm{ and } 
\frac{d\Theta_0}{\rho^2}(.,J.)=g \textrm{ on }\ker\Big(\frac{\Theta_0}{\rho^2}\Big)\subset K_1\] 
which, using the bundle isomorphism $\psi$, is actually equivalent to the following
\[(\textrm{H4'})\quad \quad \exists\;  J\in \textrm{End}(\ker\Theta_0), \; J^2=-\textrm{Id} \textrm{ and } 
d\Theta_0(.,J.)=\rho^2g|_{\ker\Theta_0}\] 
where the restriction $\rho^2g|_{\ker\Theta_0}$ 
is the metric on the bundle $\ker\Theta_0\subset T\pl\bar{X}$ 
whose value on fiber $(\ker\Theta_0)_p$ is the limit 
$\lim_{t\to 0}t^2(\phi^*g)_{(t,p)}$.

\subsection{Asymptotically complex hyperbolic manifolds}\label{ACHM}

An \emph{aymptotically complex hyperbolic manifold,} or \emph{ACH manifold,} is a non-compact Riemannian manifold 
$(X,g)$ such that there exists a smooth compact manifold with boundary $\bar{X}$ 
which compactifies $X$, equipped with a $\Theta$-structure, such that $g$ is a 
$\Theta$-metric satisfying assumptions $(\textrm{H}1)$ to $(\textrm{H}4).$

In view of the above discussion there exists a product decomposition $(0,\eps)_\rho\x
\pl\bar{X}$ near the boundary where the metric can be expressed by 
\begin{equation}\label{modform3}
g=\frac{4d\rho^2+h(\rho)}{\rho^2}
\end{equation}
with $h(\rho)$ a smooth family of metrics on $\pl\bar{X}$ for $\rho\not=0$ such that 
\begin{equation}\label{modform4}
\frac{h(\rho)}{\rho}=\frac{\Theta_0^2}{\rho^4}+\frac{d\Theta_0(.,J.)}{\rho^2}
+\rho\til{g}, \quad \til{g}\in C^{\infty}(\bar{X},S^2(^\Theta T^*\bar{X})).
\end{equation}

The form $\Theta_0$ induces the metric $h_0$ and the volume density on $\pl\bar{X}$
\begin{equation}\label{h0}
h_0:=\Theta_0^2+d\Theta_0(.,J.),\quad \textrm{dvol}_{h_0}=|\Theta_0\wedge d\Theta_0^{n}|.
\end{equation}
By choosing a different representative $\hat{\Theta}_0=e^{\omega_0}\Theta_0$,
it is easy to check that the corresponding metric is 
$\hat{h}_0=e^{2\omega_0}\Theta_0^2+e^{\omega_0}d\Theta_0(.,J.)$ and the 
volume form is $\textrm{dvol}_{\hat{h}_0}=e^{(n+1)\omega_0}\textrm{dvol}_{h_0}$.
It is then natural to call the couple $([\Theta_0],J)$ a \emph{conformal pseudohermitian structure} on $\pl\bar{X}$. 

In view of the assumptions on $\Theta_0$ for an ACH manifold, 
there exists a smooth global vector field, denoted $T_0$, tangent to $\pl\bar{X}$
such that $\Theta_0(T_0)=1$ and $d\Theta_0(T_0,Y)=0$ for every $Y\in\ker\Theta_0$; 
this is \emph{Reeb's vector field}. 
With the notation of (\ref{modform3}), we can define 
\begin{equation}\label{krho}
k(\rho):=\frac{M^{*}_\rho h(\rho)}{\rho^2},\end{equation}
where $M_\delta:T\pl\bar{X}\to T\pl\bar{X}$ is the dilation $M_\delta(tT_0+V):=\delta^2tT_0+\delta V$
if $V\in\ker\Theta_0$, $t\in\rr$ and $T_0$ is Reeb's vector field.
Observe that $k(\rho)$ is a smooth family of metrics on $\pl\bar{X}$ up to $\rho=0$, $k(0)=h_0$ and the volume form is $\textrm{dvol}_g=\rho^{-2n-2}d\rho \textrm{dvol}_{h(\rho)}=\rho^{-2n-3}d\rho \textrm{dvol}_{k(\rho)}$.\\

\textsl{Remark}: 
An ACH manifold in the sense of Biquard-Herzlich \cite{BH}
is quite similar to our setting, the difference lies in the term $\til{g}$
of (\ref{modform4}): for them, $\til{g}=O(\rho^{\delta})$ for some $\delta>0$ 
and $\til{g}$ does not have a polyhomogeneous expansion at the boundary, whereas 
in our case the metric is polyhomogeneous but we can allow terms of order 
$O(\rho^{-3})$ in the $\Theta_0$ direction, for instance. 

\section{Heisenberg pseudo-differential operators}

\subsection{$\Theta$ calculus}

We denote by $\textrm{Diff}_\Theta^m(\bar{X})$ the set of differential operators
of order $m$ which are locally polynomial functions (with coefficients 
in $C^\infty(\bar{X})$) of vector fields of $\mc{V}_\Theta$, i.e the envelopping algebra
of $\mc{V}_\Theta$. The Laplacian 
of a $\Theta$-metric is an operator in $\textrm{Diff}_\Theta^2(\bar{X})$. 
In \cite{EMM}, the authors construct a class of pseudo-differential 
operators $\Psi^*_\Theta(\bar{X})$ on $\bar{X}$ which is large enough 
to capture the resolvent $(\Delta_g-z)^{-1}$. It is defined in term 
of distributional kernel lifted on a parabolic blown-up version of $\bar{X}\x\bar{X}$\\ 

We define the blow-up (or stretched product) $\bar{X}_\Theta^2$ following Epstein-Melrose-Mendoza \cite{EMM} 
we refer the reader to sections 5,6,7 of \cite{EMM} for more details.
Let $\pi_R,\pi_L$ be the right and left projections
of $\bar{X}\x\bar{X}$ onto $\bar{X}$, and for a choice of boundary defining function $\rho$
in $\bar{X}$, we denote 
\[\rho:=\pi_L^*\rho, \quad \rho':=\pi_R^*\rho.\]
Let $\Delta$ be the diagonal in $\bar{X}\x\bar{X}$, $\pl\Delta$ its boundary, and
let $S\subset N^*(\pl\Delta)$ be the line subbundle of the conormal bundle of $\pl\Delta$ in $T^*(\bar{X}\x\bar{X})$ 
spanned by $\pi_{L}^*\Theta-\pi_{R}^*\Theta$. 
We denote by $\bar{X}^2_{\Theta}:=[\bar{X}\x\bar{X}; \pl\Delta,S]$ the $S$-parabolic blow-up of $\bar{X}\x\bar{X}$
around $\pl\Delta$. This means that we construct a larger manifold than 
$\bar{X}\x\bar{X}$ by replacing $\pl\Delta$ by the $S$-parabolically spherical normal interior pointing 
bundle in $\bar{X}\x\bar{X}$:
\[\bar{X}^2_{\Theta}:=(\bar{X}\x\bar{X}\setminus \pl\Delta)\sqcup SN_{S,+}\pl\Delta\]
where $SN_{S,+}$ means the bundle over $\pl\Delta$ such that each fiber is the quotient
of the interior pointing normal bundle (without the $0$ section) $N_+\pl\Delta\setminus\{\pl\Delta\}$ of 
$\pl\Delta\subset \bar{X}\x\bar{X}$ by the equivalence law 
\begin{equation}\label{dilationeq}
(u,z)\sim (u',z')\iff \exists \delta>0, (u,z)=M_\delta(u',z'):=(\delta u',\delta^2 z')\end{equation}
after decomposing $N_+\pl\Delta=S^0\oplus S'$ for $S^0$ annihilator of $S$ and $S'$ a 
complementary space. The complementary space has to be taken a certain way so that the total space is invariant under 
the dilation $M_\delta$ of (\ref{dilationeq}), but this can be done (see \cite[Sec. 6]{EMM}). 
We call $SN_{S,+}\pl\Delta$ the front face and denote it ff, we define the blow-down map 
\[\beta: \bar{X}^2_{\Theta}\to \bar{X}\x\bar{X}\]
to be the identity outside ff and the projection on the base on ff (recall ff is a bundle). 
We can put a topology and smooth structure on $\bar{X}^2_{\Theta}$ by taking those of $\bar{X}\x\bar{X}$ far 
from $\pl\Delta$ via $\beta$ and near ff, one can use a choice of normal fibration 
\[f: U\to V, \quad f|_{\pl\Delta}=\textrm{Id}
, \quad df=\textrm{Id} \textrm{ on } N_+\pl\Delta\]
where $U$ (resp. $V$) is a small neighbourhood of the submanifold $\pl\Delta$ in $N_+\pl\Delta$
viewed as zero section (resp. $\bar{X}\x\bar{X}$), and transport the topology and smooth structure of
$[N_+\pl\Delta; \pl\Delta,S]$ through $f$ after we have defined them on this space as follows:
homogeneous functions of non-negative integer order in $N_+\pl\Delta\setminus \pl\Delta$ with respect
to dilation $M_\delta$ lift under $\beta^*{f^{-1}}^*$ to well defined functions on $\beta^{-1}(V)$, 
the topology is the weakest such that these are continuous functions, the smooth structure
is that generated by these homogeneous functions.
This is proved in \cite{EMM} to be invariant with respect to choices of $f,S'$ as a smooth manifold
with corners.
The manifold $\bar{X}^2_{\Theta}$ has three boundary hypersurfaces, the front face ff,
the right boundary rb and the left boundary lb (cl means closure)
\[\textrm{rb}:=\textrm{cl}\Big(\beta^{-1}(\{\rho'=0\}\setminus \pl\Delta)\Big),
\quad \textrm{lb}:=\textrm{cl}\Big(\beta^{-1}(\{\rho=0\}\setminus \pl\Delta)\Big)\]
we denote by $\rho_{\frf},\rho_\lb,\rho_\rb$ some corresponding boundary defining function,
for instance one can take $\rho_\lb:=\beta^*(\rho)/\rho_\frf$, $\rho_\rb:=\beta^*(\rho')/\rho_\frf$
if $\rho_{\frf}$ is a chosen boundary defining function for ff.  
If one has a coordinate system $(\rho,u,z)$ with $(u,z)$ local coordinates on $\pl\bar{X}$
such that $\Theta_0=du+y.dx -x.dy$ is in Darboux form and $z=(x,y)\in\rr^n\x\rr^n$,
denoting the same with prime for left factor on $\bar{X}\x\bar{X}$, then 
\[\rho_\frf:=\Big((\rho^2+{\rho'}^2+\demi|z-z'|^2)^2+u^2\Big)^\frac{1}{4}, \quad \rho_{\lb}:=\frac{\rho}{\rho_\frf}, \quad
\rho_{\rb}:=\frac{\rho'}{\rho_\frf}, \quad t=\frac{u}{\rho_{\frf}^2}, \quad \omega=\frac{z-z'}{\sqrt{2}\rho_\frf}\]
are ``local coordinates'' near the front face.\\ 
 
If $\til{S}$ is the subspace of the conormal bundle $N^*\pl\Delta$ of $\pl\Delta\subset\pl\bar{X}\x\bar{X}$
spanned by $\pi_L^*\Theta_0-\pi_R^*\Theta$, then we can define similarly $\pl\bar{X}\x_\Theta\bar{X}=[\pl\bar{X}\x\bar{X};\pl\Delta,\til{S}]$ 
the $\til{S}$ parabolic blow-up of $\pl\bar{X}\x\bar{X}$ along 
$\pl\Delta$, we let $\beta'$ be the blow-down map. Finally the blow-up
$\pl\bar{X}^2_\Theta=[\pl\bar{X}\x\pl\bar{X};\pl\Delta,S_\pl]$ is defined the same way
if $S_\pl$ is the subspace of the conormal bundle $N^*\pl\Delta$ of $\pl\Delta\subset \pl\bar{X}\x\pl\bar{X}$ 
spanned by $\pi_L^*\Theta_0-\pi_R^*\Theta_0$, the blow-down map is denoted 
$\beta_\pl$.\\

Next we define the full class of $\Theta$-pseudodifferential operator
$\Psi^*_\Theta(\bar{X})$ as in Section 12 of \cite{EMM}.
An operator $A\in\rho^{E_\frf}_{\frf}\Psi_\Theta^{s,E_{\lb},E_{rb}}(\bar{X})$ 
with $(s,E_{\frf},E_\lb,E_\rb)\in \cc^{4}$ if its Schwartz kernel 
lifts under $\beta$ to a distribution $\kappa_A$  on $\bar{X}\x_\Theta\bar{X}$ such that
\[\kappa_A\in \rho_\frf^{E_{\frf}}\rho_{\lb}^{E_\lb}\rho^{E_\rb}_{\rb}C^{\infty}(\bar{X}^2_\Theta)+
\rho_{\frf}^{E_{\frf}}\mc{I}^{s}(\Delta_\iota, \bar{X}^2_\Theta)\]
where $\mc{I}^{s}(\Delta_\iota, \bar{X}\x_\Theta\bar{X})$ means the set of classically conormal distribution
of order $s$ to the interior diagonal 
\[\Delta_\iota:=\textrm{cl}\Big(\beta^{-1}(\Delta\setminus \pl\Delta)\Big)\]  
and vanishing at all order at all other boundary faces than ff. This is well-defined
since $\Delta_\iota$ meets the boundary of $\bar{X}^2_\Theta$ only at ff (in the interior
of ff). We also define $\Psi^{E,E'}(\bar{X})$ to be the set of operators
with Schwartz kernel in $\rho^{E}{\rho'}^{E'}C^{\infty}(\bar{X}\x\bar{X})$.

\subsection{The Normal operator}\label{normalsection}

If we look at the model case of $\hh^{n+1}_\cc$, it is clear that
the Lie algebra $^\Phi T_0\bar{X}$ (the boundary point is $0\in \hn$)
is canonically isomorphic to the Lie algebra of the group $\hh^{n+1}_\cc$ with the law
(\ref{grouplawhyp}), which is the reason of the notation for this last Lie algebra 
in Section 1.\\

For $p\in\pl\bar{X}$, consider $X_p:=\{d\rho\geq 0\}\subset {T_p\bar{X}}$ the inward pointing
part of $T_p\bar{X}$, then  $^\Theta T_p\bar{X}$ acts on 
$C^\infty({X_p})$ as follows.
On a neighbourhood $U_p\subset \bar{X}$ of $p$, we first define the dilation 
\[M_\delta: T\bar{X}|_{U_p}\to T\bar{X}|_{U_p}, \quad M_\delta(V):=\delta \Big(d\rho(V)\pl_\rho+\pi_{\ker\Theta_0}(V)\Big)
+\delta^2 \Theta_0(V)T_0\]
where $\pi_{\ker\Theta_0}$ is the projection on $\ker\Theta_0$ parallel to $T_0,\pl_\rho$.
Let $V(p)\in {^\Theta T}_p\bar{X}$ such that $(\rho^{-1}d\rho)(V(p))>0$, let $V\in\mc{V}_\Theta$ whose 
value at $p$ is $V(p)$, and set $V'(p)=\lim_{\delta\to 0}M_\delta^{*}V$
which is a well defined smooth vector field on $X_p$, 
homogeneous of degree $0$ with respect to $M_\delta$ and
depending only on $V(p)$.
Then the map 
\begin{equation}\label{map1}
 \left\{\begin{array}{ccc}
^\Theta T_p\bar{X}& \to &C^\infty(X_p,TX_p)\\
V(p) & \to & V'(p)
\end{array}
\right., \end{equation}
realizes $^\Theta T_p\bar{X}$ as a Lie algebra of smooth vector fields on $X_p$,
left invariant for the group action on $X_p$ generated by $^\Theta T_p\bar{X}$, 
this action on $X_p$ being 
\begin{equation}\label{groupactp}
V._pV':= d\rho(V)d\rho(V')\pl_\rho + 
\Big(\Theta_0(V+M_{d\rho(V')}V')-d\Theta_0(W,M_{d\rho(V')}W')\Big)
T_0+W+M_{d\rho(V')}W'\end{equation}
\[W=\pi_{\ker\Theta_0}(V), \quad W'=\pi_{\ker\Theta_0}(V')\]
with neutral $e_p:=\pl_\rho$. The value of the left-invariant vectors on $X_p$ at $e_p$ 
identifies $T_{e_p}X_p$ with the Lie algebra $^\Theta T_p\bar{X}$ of $X_p$ as usual. 
Using 
\[\Theta_0([Y(p),Z(p)])=-d\Theta_0(Y(p),Z(p))=\rho^2g|_{\ker\Theta_0}(Y(p),JZ(p)), \quad 
Y(p),Z(p)\in(\ker\Theta_0)_p,\]
with assumption $(\textrm{H}4')$, there exists an orthonormal  
basis $(X'_1,\dots,X'_n,Y'_1,\dots,Y'_n)$ of $(\ker\Theta_0)_p$ with respect to 
$\rho^2g|_{\ker\Theta_0}$such that 
$JX_i'=-Y_i'$.
The set $(\hh^{n+1}_\cc)_0=\{d\rho_0>0\}\subset T_0\bar{\hh}^{n+1}_\cc$ is clearly identified as a group 
with $\hh^{n+1}_\cc$ and the group action on $X_p$ is isomorphic to the law (\ref{grouplawhyp}) 
of $\hh^{n+1}_\cc$ through $A_p:X_p\to (\hh^{n+1}_\cc)_0\simeq \hh^{n+1}_\cc$ defined by 
\[A_p\Big(a\pl_\rho+bT_0+\sum_{i=1}^nx_iX_i'+y_iY_i'\Big):=a\pl_{\rho_0}+b\pl_u+\sum_{i=1}^nx_iX_i+y_iY_i\]
using notations of Section \ref{modelcase}.
This gives for each $p\in\pl\bar{X}$ an isometric Lie algebra 
isomorphism $(^\Theta T_p\bar{X},g_p) \to (^\Phi T_0\hh^{n+1}_\cc,g_{0})$ 
given by $(dA_p)_{e_p}=A_p$ (actually the linear extension of $A_p$ to $T_p\bar{X}$)  
such that $A_p^*\Phi=\Theta$.  

The normal operator $N_p$ is the map from the envelopping algebra 
$\textrm{Diff}_\Theta(\bar{X})$
of $\mc{V}_\Theta$ to the envelopping algebra 
$\mc{D}(^\Theta T_p\bar{X})$ of $^\Theta T_p\bar{X}$
induced by the projection $V\to V(p)$ from 
$\textrm{Diff}_\Theta(\bar{X})$ to $^\Theta T_p\bar{X}$. 
Then for $A\in\textrm{Diff}_\Theta(\bar{X})$, $N_p(A)$ acts by convolution 
as a left invariant differential operator on $X_p$ by (\ref{map1}).
The normal operator can be considered as an operator on $\hh^{n+1}_\cc$ by conjugating
$(A_p^{-1})^*VA_p^*$ if $V\in {^\Theta T_p\bar{X}}$.\\

The metric $g_p$ on ${^\Theta T}_p\bar{X}\simeq T_{e_p}X_p$ induces 
a left-invariant metric on $X_p$ and it is easy to check that  
this is the metric obtained by $\lim_{\delta\to 0}M_\delta^*g$.
A computation leads to 
\[N_p(\Delta_g)=\Delta_{g_p}=-\frac{1}{4}(\rho\pl_\rho)^2+\frac{n+1}{2}\rho
\pl_{\rho}-\rho^4T_0(p)^2+\rho^2\Delta_b(p)\]
where $\Delta_b(p)$ is the horizontal sublaplacian on $\pl X_p$ equipped with 
the Reeb field $T_0(p)$ and the contact distribution $(\ker\Theta_0)_p$.
Using conjugation with $A_p$, this is the complex hyperbolic Laplacian
on $\hh^{n+1}_\cc$.\\

If we look at the action of vector fields in $\mc{V}_\Theta$ lifted 
through  $\beta^*\pi_L^*$ on $\bar{X}^2_\Theta$, this action 
restricts smoothly on each fiber of the front face $\frf_p$ 
to be the left invariant action of $^\Theta T_p\bar{X}$ on $X_p$ with an identification 
between  $\frf_p$ and $X_p$. First set $\bar{X}_p:=\{d\rho\geq 0\}\subset T_p\bar{X}$, then 
\[\frf_p= \Big(((\bar{X}_p\x \bar{X}_p)/T_{p,p}\pl\Delta)\setminus\{0\}\Big)/M_\delta \]
$M_\delta$ being a parabolic dilation in $N^+_p\pl\Delta$ induced 
by the kernel of $\pi_R^*\Theta-\pi_L^*\Theta$ and a choice of transversal, 
here $\pi_R^*T_0-\pi_L^*T_0$.
If $s,u,v \to s\pl_\rho+uT_0+v$ with $\Theta_0(v)=0$ are coordinates on $\bar{X}_p$
then the change of coordinates $(s,u,v;s',u',v')\to (s,t=u-u',z=v-v';s',u',v')$ gives an isomorphism 
of $X_p\x X_p$ and we have $T_{p,p}\pl\Delta=\{s=t=z=s'=0\}$ thus
\[\frf_p\simeq \Big((\bar{X}_p\x [0,\infty))\setminus\{0\}\Big)/M'_\delta, \quad M'_\delta(s,t,z; s'):=(\delta s,\delta^2t,\delta z; \delta s')\]
with $\{0\}=\{s=t=z=s'=0\}$. Setting $|z|^2:=d\Theta_0(z,Jz)$, the identification of $\frf_p$ with the quarter of sphere 
$\rho_\frf:=((s^2+{s'^2}+|z|^2/2)^2+t^2)^{\frac{1}{4}}=1$ gives the smooth structure of manifold with corners and 
the functions $\rho_\lb:=s/\rho_\frf$, $\rho_\rb:=s'/\rho_\frf$ are boundary defining functions for 
left and right faces. If we fix $s'=1$, we have an identification between $X_p$ and the 
interior of $\frf_p$ given by ``stereographic projection'' with pole $\rho_\rb=1$ 
for the parabolic dilation, these are projective coordinates of $\frf_p$. We then 
have in this model of $\frf_p$
\begin{equation}\label{rholbrb}
\rho_\lb=\frac{s}{((s^2+1+\demi|z|^2)^2+t^2)^{\frac{1}{4}}}, \quad 
\rho_\rb=\frac{1}{((s^2+1+\demi|z|^2)^2+t^2)^{\frac{1}{4}}}.
\end{equation}
A vector in $V_p\in {^\Theta T}_p\bar{X}$, value at $p$ of a vector field $V\in\mc{V}_\Theta$, 
acts as a homogeneous (for $M_\delta$) left invariant vector field 
on $X_p$, thus on $\frf_p$, by the dilated limit (\ref{map1}). 
This left invariant field acting on $X_p$ 
is $N_p(V)$ but this is also the restriction of $\beta^*\pi_R^*V$ to the interior of $\frf_p$.  
It is proved in \cite{EMM} that a vector field $V\in\mc{V}_\Theta$
actually lifts smoothly to $\bar{X}^2_\Theta$ under $\beta^*\pi_R^*$, tangent to any 
boundary hypersurface.

\subsection{Heisenberg pseudo-differential operators on Heisenberg manifolds}\label{heisenpseudo}

The space of Heisenberg pseudo-differential operators is 
defined in \cite{BG} (see also the monograph of Ponge \cite{Po}) 
but the approach we will use is that of Epstein-Melrose-Mendoza \cite{EMM2} since it is more naturally 
adapted to our case. In \cite{EMM2} they define the class $\Psi^m_{\Theta_0}(\pl\bar{X})$
of classical Heisenberg pseudodifferential operators of order $m$ 
by the structure of their Schwartz kernel.\\ 

Let $M$ be a compact manifold of dimension $n$
equipped with a one-dimensional subbundle $L=\rr\Theta_0\subset T^*M$. 
one defines ${^L\bar{T}}^*M$ to be the 
parabolic compactification of $T^*M$ with respect to $L$, this is a smooth manifold 
with boundary and let $q$ be 
a boundary defining function. A Heisenberg classical symbol of order $m$ 
is a function in $q^{-m}C^\infty({^L\bar{T}}^*M)$.  
If $\Delta$ is the diagonal of $M\x M$, we denote by $N\Delta$ the normal 
bundle of $\Delta$ in $M\x M$ and there is a canonical Fourier transform
of compactly supported distributions in each fiber and a neihgbouhood of 
the diagonal in $M\x M$ can be identified to a neighbourhood of the zero section 
of $N\Delta$. Note also that $I:V\to (\demi V,-\demi V)$ identifies $T_xM$ with a subspace 
of $T_xM\oplus T_xM$ cannonicaly isomorphic with $N_x\Delta$. 
A choice of metric $h$ yields the Riemann-Weyl fibration $W:U\subset TM\to M\x M$ 
\[W: (x,V)\to (\exp_x^h(\demi V),\exp_{x}^h(-\demi V))\]
where $U$ is a neighbourhood of the zero section, we define by $(m,\tau)=W^{-1}$ and
$\psi:=W\circ I^{-1}$ normal fibration $I(U)\subset N\Delta\to M\x M$.
An operator $A$ is in $\Psi^m_{\Theta_0}(M)$ if its Schwartz kernel
$k$ is smooth outside the diagonal and can be written near the diagonal under the form
\begin{equation}\label{kxy}
k(x,y)=(2\pi)^{-n}\int_{T^*_{m(x,y)}}e^{i\tau(x,y).\xi}\sigma(m(x,y),\xi)d\xi
,\quad \sigma\in \rho^{-m}C^\infty({^L\bar{T}}^*M).
\end{equation}
The principal symbol is the class $\sigma_0:=\in q^{-m}C^\infty({^L\bar{T}}^*M)/
q^{-m+1}C^\infty({^L\bar{T}}^*M)$ of $\sigma$. Epstein-Melrose-Mendoza \cite{EMM2} proves that 
both definitions do not depend on the choice of $h$. 

We can now define $S$ to be the span of $\pi_R^*\Theta_0-\pi_L^*\Theta_0$ in 
$N^*\Delta$. We let $M^2_{\Theta_0}:=[M\x M; \Delta,S]$ be the parabolic 
blow-up of $M\x M$ around $\Delta$ in direction $S$, $\rho_\frf$ 
be a boundary defining function of the front face $\frf$ and $\beta$ the blow-down map.
Let $Y:=[N\Delta; \Delta, S]$ and $\beta_Y$ the associated blow-down map, 
$Y$ and $M^2_{\Theta_0}$ have same front face and $\beta=\beta_Y$ on this face.  
Moreover $\psi=\textrm{Id}$ on $\Delta$, $\psi_*=\textrm{Id}$ on $N\Delta$
and $\psi$ lifts smoothly (by taking the identity on $\frf$) as a local 
diffeomorphism $\til{\psi}: Y\to M^2_{\Theta_0}$ near $\Delta$ by construction of 
the blow-up smooth structure. Moreover the way that $\psi$ is constructed 
(by Riemmann-Weyl fibration) implies the important identity
\begin{equation}\label{dpsi}
d\psi=\textrm{Id} \textrm{ on } \frf.
\end{equation}
Let $\rho_{\frf(Y)}$, $\rho_\frf$ be boundary defining function for respectively 
$\frf\subset Y$ and $\frf\subset M^2_{\Theta_0}$, 
then for any $s\in\cc$, functions in $\rho_{\frf}^sC^\infty(M^2_{\Theta_0})$ 
defined near $\frf$ lift under $\til{\psi}$ to functions in 
$\rho_{\frf(Y)}^sC^\infty(Y)$ and the spaces 
\[\rho_{\frf(Y)}^sC^\infty(Y)/\rho_{\frf(Y)}^{s+1}C^\infty(Y),\quad  
\rho_{\frf}^sC^\infty(M^2_{\Theta_0})/\rho_{\frf}^{s+1}C^\infty(M^2_{\Theta_0})\] 
do not depend on $\rho_{\frf(Y)},\rho_\frf$ and are in one to one correspondance 
(these are conormal densities to $\frf$) which does not depend
on the choice of Riemann-Weyl fibration in view of (\ref{dpsi}).
Now if $\til{K}\in \rho_{\frf(Y)}^sC^\infty(Y)$ for $s\in\cc\setminus(-n-1-\nn)$,
then $k:={{\beta_Y}_*}\til{K}$ has an expansion near $0$ in each fibre in 
$k\sim \sum_ik_i$ with $k_i$ homogeneous of order $s+i$ 
with respect to the parabolic dilation in the fiber 
and can thus be written like (\ref{kxy}) 
for $\sigma\in q^{-s-n-1}C^\infty({^L\bar{T}}^*M)$ 
the inverse Fourier transform of ${\beta_Y}_*\til{K}$ in the fibres (Fourier 
transform keeps parabolic homogeneity).
Choosing $\rho_{\frf(Y)}$ homogeneous of order $1$ in the fibre, 
a term $\til{K}\in\rho_{\frf(Y)}^sC^\infty(Y)/\rho_{\frf(Y)}^{s+1}C^\infty(Y)$
can be uniquely represented by a homogeneous distribution of order $s$
with respect to the parabolic dilation in the fibre by taking 
\begin{equation}\label{restric} 
\rho_{\frf(Y)}^{s}(\rho_{\frf(Y)}^{-s}\til{K})|_{\frf},
\end{equation} 
its Fourier transform thus defines 
an element $\sigma_0\in q^{-s-n-1}C^\infty({^L\bar{T}}^*M)
/q^{-s-n}C^\infty({^L\bar{T}}^*M)$
independent on the choice of homogeneous $\rho_{\frf(Y)}$ and Riemann-Weyl fibration, 
i.e. the principal symbol. 
But if $K\in\rho^s C^\infty(M^2_{\Theta_0})$ then (\ref{dpsi}) insures that
(\ref{restric}) is also equal to $\rho_{\frf(Y)}^s(\rho_{\frf}^{-s}K)|_{\frf}$ 
(understood in normal coordinates $\frf\x[0,\infty)_{\rho_{\frf(Y)}}$) 
if $\til{K}:=\til{\psi}^*K$ and $\rho_\frf:=\til{\psi}_*\rho_{\frf(Y)}$, 
then the principal symbol of the operator is
\[ \sigma_0(x,.) := \mc{F}_{N_x\Delta}{\beta_Y}_*\Big(\rho_{\frf(Y)}^{s}
(\rho_\frf^{-s} K)|_{\frf_x}\Big)
=\mc{F}_{T_xM}I^*{\beta_Y}_*\Big(\rho_{\frf(Y)}^{s}
(\rho_\frf^{-s} K)|_{\frf_x}\Big)\]
where $\mc{F}$ means Fourier transform. 

Notice that the Fourier transform on tangent and cotangent bundle 
needs to be given a density to be well defined, this is 
for instance the case if one has a metric on $M$.
Note also that the definitions of $\Psi^*_{\Theta_0}(M)$ being local, the 
metric $h$ for the Riemann-Weyl fibration can be locally taken to be flat.  Finally we remark that the
contact form $\Theta_0$ and the almost-complex structure $J$ on its kernel 
induce a metric on $M$ by (\ref{h0} and a Reeb vector field $T_0$ fixing 
a representative for the parabolic blow-ups.

\section{The Laplacian, its resolvent and the Poisson operator}

We use the form (\ref{modform3}) of the metric in the product decomposition $[0,\eps)_\rho\x \pl\bar{X}$ 
of $\pl\bar{X}$ and the dilated metric $k(\rho)$ defined in (\ref{krho}).
We first write the Laplacian near infinity using this model form of the metric  
\begin{equation}\label{formedelta}
\Delta_g=-\frac{1}{4}(\rho\pl_\rho)^2+\frac{(n+1)}{2}\rho\pl_\rho+\rho^2\Delta_{h(\rho)}+
-\frac{1}{8}\rho\pl_\rho(\log|k|)\rho\pl_\rho
\end{equation}
where $|k|:=|\det k(\rho)|$, this gives 
\begin{equation}\label{decompdelta}
\Delta_g=-\frac{1}{4}(\rho\pl_\rho)^2+\frac{(n+1)}{2}\rho\pl_\rho-\rho^4T_0^2+
\rho^2\Delta_b + \rho P,\quad 
P \in\textrm{Diff}_\Theta^2(\bar{X})
\end{equation}
where $\Delta_b$ is the horizontal sub-laplacian on $(\pl\bar{X},\Theta_0,h_0)$.\\ 

\textsl{Remark}: If the metric is even at order $2k$ as defined in Subsection \ref{modelbdfs}, 
then $\nabla^{h(\rho)}f$ is even in $\rho$ at order $2k$ as well as the divergence $\textrm{div}_{h(\rho)}X=\textrm{div}_{k(\rho)}X$ 
for any $X\in T\pl\bar{X}$. Finally we have $|\det k(\rho)|=\rho^{2}|\det h(\rho)|$
which is then even in $\rho$ at order $2k+2$, this clearly proves that (\ref{formedelta}) is (in local coordinates near $\pl\bar{X}$) 
an operator of the form $Q(\rho,y;\rho\pl\rho,\pl_y)$ where $Q(x,y;D)$ is a family of polynomial in $D$
with an even Taylor expansion at $x=0$ at order $2k$.

\subsection{The Resolvent}\label{theres}

The Laplacian $\Delta_g$ has essential spectrum 
$\sigma_{\textrm{ess}}(\Delta_g)=[\frac{(n+1)^2}{4},\infty)$
and possibly a finite set $\sigma_{pp}(\Delta_g)$ of eigenvalues 
in $(0,\frac{(n+1)^2}{4})$. 
Let us set for $k\in\nn_0$
\[\mc{P}_k:=\frac{n+1/2}{2}-k-\demi\nn_0,\]
then Epstein, Melrose and Mendoza \cite{EMM} proved the 
\begin{theo}\label{emm}
On an asymptotically complex hyperbolic manifold whose metric is even at order $2k$, 
the modified resolvent of the Laplacian
\[R(\la)=(\Delta_g-\la(n+1-\la))^{-1}\]
extends from $\{\Re(\la)>\frac{n+1}{2}\}$ to $\cc\setminus (-\nn_0\cup\mc{P}_k)$ 
as a finite-meromorphic family of operators in the $\Theta$-calculus: 
for $\la$ not a pole, we have
\[R(\la)\in \Psi_{\Theta}^{-2,2\la,2\la}(\bar{X}),\quad
R(\la): \dot{C}^\infty(\bar{X})\to \rho^{2\la}C^\infty(\bar{X})\] 
\end{theo}
The poles are called resonances and we denote by $\mc{R}$ the set of resonances.\\
 
\textsl{Remark}: in \cite{EMM},  the property of metrics even at order
$2k$ is only discussed for $k=\infty$, but this can be checked in general from 
their construction (see also \cite{Gth} for a similar analysis in the real asymptotically hyperbolic case).\\
  
We define the multiplicity of a resonance $\la_0$ by 
\begin{equation}\label{defmult}
m(\la_0):=\rang \Big(\textrm{Res}_{\la=\la_0}((2\la-n-1)R(\la))\Big)
\end{equation}  
which in turn is equivalent to the rank of the polar part of $R(\la)$ at $\la_0$. 
  
\subsection{The Poisson Operator} \label{poissonsection}
   
Using the form (\ref{formedelta}) of the Laplacian and the structure of the resolvent, we are able
to define the Poisson operator:
\begin{prop}\label{operpoisson}
Suppose the metric is even at order $2k$ and let $\la$ satisfy $\Re(\la)\geq \frac{n+1}{2}$, $\la\notin 
(\frac{n+1}{2}+k+\frac{1}{4}\nn)$ and $\la(n+1-\la)\notin\sigma_{pp}(\Delta_g)$. 
Then for each $f_0\in C^\infty(\pl\bar{X})$ there exists a unique function 
$F(\la,f_0)$, linear in $f_0$, 
solving the problem
\begin{equation}\label{poissonprob}
\left\{\begin{array}{l}
(\Delta_g-\la(n+1-\la))F(\la,f_0)=0\\
F(\la,f_0)=\rho^{2(n+1-\la)}F_1(\la,f_0)+\rho^{2\la}F_2(\la,f_0)\\
F_i(\la,f_0)\in C^\infty(\bar{X})\\
F_1(\la,f_0)|_{\pl\bar{X}}=f_0
\end{array}\right.\end{equation}
with $F_{1}(\la,f_0)$ having an even Taylor expansion at $\rho=0$ at order $2k$ in a product 
decomposition near $\pl\bar{X}$.
The Poisson operator defined by 
\[P(\la): \left\{
\begin{array}{ccc}
C^\infty(\pl\bar{X}) & \to & C^\infty(X)\\
f_0 & \to  & F(\la,f_0)\end{array}\right..\]
extends finite-meromorphically to $\cc\setminus(-\nn_0\cup\mc{P}_k)$. 
\end{prop} 
\textsl{Proof}: The proof is similar to the proof of Graham-Zworski \cite{GRZ}
in the asymptotically hyperbolic case. We first deduce from  (\ref{decompdelta})
the indicial equation for $f_0\in C^\infty(M)$ and $j\in\nn_0$
\begin{equation}
(\Delta_g-\la(n+1-\la))\frac{2\rho^{2(n+1-\la)+j}f_0(y)}{j(2\la-n-1-j/2)}=\rho^{2(n+1-\la)+j}
f_0(y)+O(\rho^{2(n+1-\Re(\la))+j+1})
\end{equation}
Then we can construct a meromorphic function $\Phi(\la)\in \rho^{2(n+1-\la)}C^\infty(\bar{X})$ with at most 
first order poles at $\frac{n+1}{2}+\frac{1}{4}\nn$ such that
\[(\Delta_g-\la(n+1-\la))\Phi(\la,f_0)=O(\rho^\infty), \quad \rho^{-2(n+1-\la)}
\Phi(\la,f_0)|_{\pl\bar{X}}=f_0.\]
We thus define 
\begin{equation}\label{constpoisson}
P(\la)f_0=\Phi(\la,f_0)-R(\la)(\Delta_g-\la(n+1-\la))\Phi(\la,f_0)
\end{equation}
which satisfies the required conditions. However we have to prove uniqueness 
of an operator satisfying (\ref{poissonprob}), this is obtained easily using the 
indicial equation when $\Re(\la)>\frac{n+1}{2}$ and 
using a Green formula at infinity if $\Re(\la)=\frac{n+1}{2}$ (see \cite{GRZ} for 
details) :
\begin{lem}\label{boundarypairing}
Let $\Re(\la)=\frac{n+1}{2}$ and $u_i=\rho^{2(n+1-\la)}F_i+\rho^{2\la}G_i$ for some 
$F_i,G_i\in C^\infty(\bar{X})$ and $i=1,2$. If $u_i$ satisfies
$r_i:=(\Delta_g-\la(n+1-\la))u_i \in \dot{C}^\infty(\bar{X})$
for $i=1,2$, then 
\[\int_X (u_1\bbar{r_2}-u_2\bbar{r_1})\textrm{dvol}_{g}=2(2\la-n-1)\int_{\pl\bar{X}}
(F_1|_{\pl\bar{X}}\bbar{F_2}|_{\pl\bar{X}}-G_1|_{\pl\bar{X}}\bbar{G_2}|
_{\pl\bar{X}})\textrm{dvol}_{h_0}.\] 
\end{lem}
\textsl{Proof}: we write Green's formula in $\{\rho>\eps\}$ 
\[\int_{\rho>\eps}(u_1\bbar{r_2}-u_2\bbar{r_1})\textrm{dvol}_{g}\eps^{-2n+1}\int_{\rho=\eps}(u_1\pl_\rho\bbar{u_2}-\bbar{u_2}\pl_\rho u_1)\textrm{
dvol}_{h(\rho)}.\]
Using $\textrm{dvol}_{h(\rho)}=\rho^{-1}\textrm{dvol}_{h_0}+O(1)$ and the asymptotics 
of $u_i$ gives
\[u_1\pl_\rho\bbar{u_2}-\bbar{u_2}\pl_\rho u_1=2(2\la-n-1)\rho^{2n}(F_1|_{\pl\bar{X}}
\bbar{F_2}|_{\pl\bar{X}}
-G_1|_{\pl\bar{X}}\bbar{G_2}|_{\pl\bar{X}})+O(\rho^{2n+1})\]
which implies the Lemma by taking the limit as $\eps\to 0$.
\qed\\

Then the first part of the proposition is deduced classically from this Lemma and 
formula (\ref{constpoisson}) extends clearly meromorphically to 
$\cc\setminus(-\nn_0\cup\mc{P}_k\cup (n+1-\mc{P}_0))$
using Theorem \ref{emm}, which proves the meromorphic extension of $P(\la)$.
The points $(n+1/2)/2+j/2$ are not singularities if $j<2k$ since by evenness of the metric 
at order $2k$ and the remark following (\ref{decompdelta}), the function $\rho^{-2(n+1-\la)}\Phi(\la,f_0)$ can be constructed 
to have even Taylor expansion at order $2k$ at $\rho=0$, see \cite{G0} where it is done in greater details
for the real asymptotically hyperbolic case.
 
>From the proof of proposition 3.5 in \cite{GRZ}, we also deduce that 
$P(\la):C^\infty(\pl\bar{X})\to C^\infty(X)$, 
is analytic at each 
\begin{equation}\label{lak}
\la_k:=\frac{n+1}{2}+\frac{1}{4}k,\quad k\in\nn 
\end{equation}
if $\la_k(n+1-\la_k)\not\in\sigma_{pp}(\Delta_g)$, and $\mc{P}(\la_k)f_0$
can be defined as the unique solution of 
\begin{equation}\label{poissonprobk}
\left\{\begin{array}{l}
(\Delta_g-\la_k(n+1-\la_k))\mc{P}(\la_k)f_0=0\\
\mc{P}(\la_k)f_0=\rho^{2(n+1-\la_k)}F_k(f_0)+\rho^{2\la_k}\log(\rho)G_k(f_0)\\
F_k(f_0),G_k(f_0)\in C^\infty(\bar{X})\\
F_k(f_0)|_{\pl\bar{X}}=f_0
\end{array}\right..\end{equation}
and we are done.
\qed\\

Actually the Poisson operator was constructed in 
Section 15 of \cite{EMM} by its Schwartz kernel as 
a weighted restriction of the resolvent kernel at
the boundary face $\textrm{lb}(\bar{X}^2_\Theta)$. 
The Schwartz kernel of $P(\la)$ is 
\begin{equation}\label{noyaup1}
P(\la)=c(\la)\beta'_*\Big(\beta^*({\rho'}^{-2\la}R(\la))
|_{\rb}\Big)\end{equation}
for some smooth positive function $c(\la)\in C^\infty(\pl\bar{X}\x\bar{X})$.
Actually one can use Lemma \ref{boundarypairing} as in \cite[Prop. 3.9]{GRZ} to check 
the same result, which gives the value $c(\la)=2(2\la-n-1)$. 
Using Theorem \ref{emm}, we deduce that  
\begin{equation}\label{noyaup2}
\beta^*({\rho'}^{-2\la}R(\la))|_{\rb}\in \rho_{\textrm{lb}}^{2\la}\rho_{\textrm{ff}}^{-2\la}
C^\infty(\bar{X}\x_\Theta \pl\bar{X})+\beta^*(\rho^{2\la}C^{\infty}(\bar{X}\x\pl\bar{X})).
\end{equation}
 
A property that we will use later is the relation between spectral measure and Poisson operator, 
the proof of which is an easy appalication of Green formula, and essentially similar to the real asymptotically hyperbolic case (see \cite{P,Gu1,G0}):
\begin{lem}\label{resvspoisson}
If the metric is even at order $2k$, we have for $\la,n-\la\notin \mc{R}\cup \mc{P}_k$
\[R(\la;m,m')-R(n+1-\la;m,m')=(2(2\la-n-1))^{-1}\int_{\pl\bar{X}}P(\la;m,y)P(n-\la;m',y){\rm dvol}_{h_0}(y)\] 
\end{lem}
 
\section{The Scattering operator}
>From the Poisson problem, we can define the scattering operator 
for $\Re(\la)\geq \frac{n+1}{2}$ and $\la\not\in
(\mc{R}\cup\frac{n+1}{2}+\frac{1}{4}\nn)$ by 
\begin{equation}\label{defscat}
S(\la):\left\{\begin{array}{ccc}
C^\infty(\pl\bar{X}) & \to & C^\infty(\pl\bar{X})\\
f_0 & \to & F_2(\la,f_0)|_{\pl\bar{X}} 
\end{array}\right.
\end{equation}
where $F_2(\la,f)$ is defined by (\ref{poissonprob}). 
Since $P(\la)$ extends finite meromorphically to $\cc\setminus (-\nn_0\cup\mc{P}_k)$ if the metric is even
at order $2k$ by Proposition (\ref{operpoisson}), we deduce that 
$S(\la)$ continues meromorphically in $\cc\setminus (-\nn_0\cup\mc{P}_k)$ 
\[S(\la)f=-\Big(\rho^{-2\la}R(\la)(\Delta_g-\la(n+1-\la))\Phi(\la,f)\Big)|_{\rho=0}.\]
The rank at a pole $\la=s$ is finite if $s\not\in -\nn_0\cup\mc{P}_k\cup(\frac{n+1}{2}+\frac{1}{4}\nn)$ 
and in $\{\Re(\la)>\frac{n+1}{2}\}$, the same arguments than \cite[Prop. 3.6]{GRZ} show that 
$S(\la)$ has at most first order poles with residue 
\[\textrm{Res}_{s}S(\la)=\left\{\begin{array}{ll}
\Pi_{s} & \textrm{ if } s\not\in\frac{n+1}{2}+\frac{1}{4}\nn\\
\Pi_{\la_k}-p_k & \textrm{ if } s=\la_k\end{array}\right.\]
with $\Pi_{s}$ a finite rank operator with Schwartz kernel
\begin{equation}\label{pis}
\pi_{s}=2(2s-n-1)
\Big((\rho\rho')^{-2s}\textrm{Res}_{s}R(\la)\Big)|_{\rho=\rho'=0}
\end{equation}
and $p_k$ is a differential operator on $\pl\bar{X}$ defined by
\begin{equation}\label{defpk}
p_j=\textrm{Res}_{\la_k}p_j(\la), \textrm{ where } 
\Phi(\la,f_0)\sim\sum_{j}\rho^{2(n+1-\la)+j}p_j(\la)f_0.\end{equation}
It is straightforward to see from the construction of $\Phi(\la)$ and (\ref{decompdelta}) 
that $p_j(\la)$ is a differential operator on $\pl\bar{X}$,
actually a Heisenberg differential operator,
the Heisenberg principal symbol $\sigma_{\rm pr}(p_{2k}(\la))$ of $p_{2k}(\la)$ 
can be obtained by the induction formula
\[p_0(\la)=1 , \quad p_2(\la)=-\Delta_b, \quad p_{2k+2}(\la)=\frac{-\Delta_b p_{2k}(\la)+
T_0^2p_{2k-2}(\la)}{(k+1)(2\la-n-1-k-1)}\textrm{ mod }\Psi_{\Theta_0}^{2k-1}(\pl\bar{X}).\]
Since $[\Delta_b,T_0]\in \Psi_{\Theta_0}^{3}(\pl\bar{X})$ and we are interested only in the principal
symbol, it suffices to deal as if $\Delta_b,T_0$ were complex numbers. 
Observing that $p_{2k}(\la)$ is holomorphic at $(n+2+k+\nn_0)/2$, one can consider 
\[q_{2k}(\la):=\frac{\Gamma(2\la-n-1)k!}{\Gamma(2\la-n-1-k)}p_{2k}(\la), \quad q_0(\la):=1, \quad q_{2}(\la):=-\Delta_b\] 
which then satisfies
\[q_{2k+2}(\la)=-\Delta_bq_{2k}(\la)+k(2\la-n-1-k)T_0^2q_{2k-2}(\la)\]
and such that $p_{2k}=((k-1)!k!2)^{-1}q_{2k}(\la_{2k})$. But this last
term $q_{2k}(\la_{2k})$ is computed by Graham \cite[Sec. 1]{GRA} and is equal to
\[q_{2k}(\la_{2k})=\prod_{l=1}^k(-\Delta_b+i(k+1-2l)T_0).\]
We have thus proved
\begin{prop}\label{principal}
The scattering operator $S(\la)$ is a finite meromorphic family of operators
on $C^{\infty}(\pl\bar{X})$ in $\cc\setminus(-\nn_0\cup\mc{P}_k\cup(n+1-\mc{P}_k)\cup(\frac{n+1}{2}+\frac{1}{2}\nn))$
and meromorphic in $\cc\setminus(-\nn_0\cup\mc{P}_k)$. 
In $\{\Re(\la)>\frac{n+1}{2}\}$, $S(\la)$ has at most first order poles with residue 
\[\textrm{Res}_{s}S(\la)=\left\{\begin{array}{ll}
\Pi_{s} & \textrm{ if } s\not\in\frac{n+1}{2}+\frac{1}{4}\nn\\
\Pi_{\la_k}-p_k & \textrm{ if } s=\la_k\end{array}\right.\]
with $\Pi_{s}$ is the finite rank operator defined in (\ref{pis}) and 
$p_k$ is the  Heisenberg differential operator on $\pl\bar{X}$ defined in (\ref{defpk})
and such that
\begin{equation}\label{prsymbp2k}
p_{2k}=\frac{1}{2((k-1)!k!)}\prod_{l=1}^k(-\Delta_b+i(k+1-2l)T_0) \textrm{ mod }\Psi_{\Theta_0}^{2k-1}(\pl\bar{X}).\end{equation}
\end{prop}

Note that $S(\la)$ depends on the boundary defining function 
$\rho$ and a different choice $\hat{\rho}=e^{\omega}\rho$ (i.e. a 
different conformal representative 
contact form $e^{2\omega}\Theta_0$ with $\omega_0=\omega|_{\pl\bar{X}}$)
gives a scattering operator $\hat{S}(\la)$ and residues $\hat{p}_k$ (when $\la(n+1-\la)\notin \sigma_{\textrm{pp}(\Delta_g)}$) conformally related to $S(\la),p_k$ 
\[\hat{S}(\la)=e^{-2\la\omega_0}S(\la)e^{2(n+1-\la)\omega_0}, \quad \hat{p}_k=e^{-(n+1+k/2)\omega_0}p_ke^{(n+1-k/2)\omega_0}.\]
The definition of $S(\la)$ and Lemma \ref{boundarypairing} directly imply (see \cite[Sec. 3]{GRZ}) that 
\begin{equation}\label{eqfct}
S(\la)^{-1}=S(n-\la)=S^{*}(\la), \quad \Re(\la)=\frac{n+1}{2} 
\end{equation}
where the adjoint is taken with respect to measure $\rm{dvol}_{h_0}$ on $\pl\bar{X}$ and the left identity
extends meromorphically in $\cc\setminus (\demi\zz\cup\mc{P}_k\cup(n+1-\mc{P}_k))$.\\

We have a direct relation between resolvent and scattering operator:
\begin{prop}\label{kernelsla}
For $\la$ such that $S(\la)$ exists, its Schwartz kernel is
\[S(\la)=2(2\la-n-1){\beta_\pl}_*\left(\beta^*((\rho\rho')^{-2\la}R(\la))|_{{\rm lb}\cap{\rm rb}}\right).\]
The right-hand side is a meromorphic family of distributions on $\pl\bar{X}$ in $\cc\setminus -\nn_0\cup\mc{P}_k$
with first order poles at $(n+1+\nn)/2\cup (n+1-\mc{P}_k)$ if the metric is even at order $2k$.
\end{prop}
\textsl{Proof}: 
For $\Re(\la)<0$, we clearly have that 
\[S(\la)f=\lim_{\rho\to 0}[\rho^{-2\la}(P(\la)f-\Phi(\la,f)]\lim_{\rho\to 0}[\rho^{-2\la}P(\la)f].\]
Using (\ref{noyaup1}) and (\ref{noyaup2}) for $\Re(\la)<0$ gives  
\[\rho^{-2\la}P(\la) \in \beta_*(\rho_{\frf}^{-4\la}C^\infty(\bar{X}\x_
\Theta\pl\bar{X}))\subset C^\infty(\bar{X}
\x\pl\bar{X}\setminus \pl\Delta)\cap C^0(\bar{X}\x\pl\bar{X})\]
thus by taking the limit 
\[S(\la)f(b)=\int_{\pl\bar{X}}\lim_{m\to b}(\rho(m)^{-2\la}\mc{P}(\la;m,b'))f(b')
\textrm{ dvol}_{h_0}(b')\] 
we deduce that the kernel of $S(\la)$ is 
${\beta_\pl}_*({\beta'}^*(\rho^{-2\la}P(\la))|_{\lb})$ for $\Re(\la)<0$,
this gives the result in a half plane by using again (\ref{noyaup1}).
To extend it in the complex plane, we need to check that
a distribution of the form $\rho_{\frf}^{-4\la}F_\la$ with 
$F_\la\in C^\infty(\pl\bar{X}\x_\Theta\pl\bar{X})$ meromorphic 
extends meromorphically to $\cc$. But this is a classical fact which 
can be proved as in the normal blow-up case (classical pseudo-differential 
operators) by Taylor expansion at the front face of the parabolic 
blow-up, we lead the reader to the Appendix where we describe it. 
This leads to first order poles at $\frac{n+1}{2}+\frac{1}{4}\nn$
whose residue are differential operators. 
\qed\\

For $\Theta_0$ fixed, we have an induced metric $h_0$ on $\pl\bar{X}$ (see (\ref{h0})), 
thus orthogonal decomposition $T\pl\bar{X}=\rr T_0\oplus \ker\Theta_0$ and 
$T^*\pl\bar{X}=\rr \Theta_0\oplus (\rr\Theta_0)^\perp$. 
We define the  Heisenberg norm on $T\pl\bar{X}$  
\[||\xi||_{\he}:=\Big(\Theta_0(\xi)^2+\frac{1}{4}d\Theta_0(\xi,J\xi)^2\Big)^{\frac{1}{4}}.\]
The Fourier transform of $||\xi||^{s}$ is $s\in\cc$ 
is a meromorphic family, with poles at $\nn$, of homogeneous distributions 
of order $-s-2n-2$ on $T^*\pl\bar{X} 
=\rr \Theta_0\oplus (\rr\Theta_0)^\perp$ 
with respect to the parabolic dilation $M_\delta(t\Theta_0+ u):=\delta^2t\Theta_0+\delta u$, with $u\perp \Theta_0 $. 

\begin{theo}\label{pseudo}
The scattering operator $S(\la)$ is, for $\la\not\in \mc{R}\cup \frac{1}{4}\zz$,
a pseudo-differential operator in the Heisenberg class 
\[S(\la)\in \Psi_{\Theta_0}^{4\la-2(n+1)}(\pl\bar{X}).\] 
and the principal symbol of $S(\la)$ is 
\[\sigma_{\rm pr}(S(\la))(\xi)=c_n\frac{2^{2\la+1}\Gamma(\la)^2}{\Gamma(2\la-n-1)}
\mc{F}_{V\to\xi}(||V||_{\rm{He}}^{-4\la})\]
\end{theo}
\textsl{Proof}: we use \cite{EMM2} definition of $\Psi^*_{\Theta_0}(\pl\bar{X})$
(that we recalled in Subsection \ref{heisenpseudo}) and the first property is a consequence of Proposition \ref{kernelsla} 
and discussions in Subsection \ref{heisenpseudo}. 

To compute the principal symbol, we work locally near a point $(p,p)$
of the diagonal $\pl\Delta$ in the boundary and follow Subsection \ref{heisenpseudo}; 
we set $M=\pl\bar{X}$ and define the bundle 
$Y:=[N\pl\Delta(M^2);\pl\Delta, S_\pl]$ where 
$S_\pl:=\pi_R^*\Theta_0-\pi_L^*\Theta_0$, and suppose $\rho_{\frf(Y)}$ is a boundary 
defining function of the front face $\frf(Y)$ of $Y$ which is homogeneous with respect to the 
parabolic dilation induced by $S_\pl$ and the Reeb field $T_0$. 
We also use $M^2_{\Theta_0}:=[M\x M;\pl\Delta, S_\pl]$, denote by $\pl\frf$ its front face 
(remark that $\pl\frf=\frf(Y)$) and denote by $\beta_\pl$ the blow-down map. 
Consider also a Riemann-Weyl fibration $\psi:N\Delta\to M^2$, 
near $\Delta$ which lifts to $\til{\psi}:Y \to M^2_{\Theta_0}$. Finally we set 
$\rho_{\pl\frf}:=\til{\psi}_*\rho_{\frf(Y)}$.
The main singularity $s_0(\la)$ of $S(\la)$ is given by the homogeneous
distribution
\begin{equation}\label{s0}
s_0(\la)=\rho_{\frf(Y)}^{-4\la}\Big(\beta_{\pl}^*(\rho_{\pl\frf}^{4\la}S(\la))|_{\pl\frf}\Big)=2(2\la-n-1)\rho_{\frf(Y)}^{-4\la}\Big([\rho_{\pl\frf}^{4\la}\beta^*((\rho\rho')^{-2\la}R(\la))|_{\lb\cap\rb}
]|_{\pl\frf}\Big)\end{equation}
and the symbol is the Fourier transform of ${\beta_\pl}_*(s_0(\la))$
in the fibers of $N\pl\Delta \simeq T\pl\bar{X}$. 
Recall that the manifold $M^2_{\Theta_0}$ is canonically identified with 
$(\bar{X}\x_{\Theta_0}\bar{X})\cap \rb\cap \lb$.
Now suppose that there exists a boundary defining function $\rho_\frf$ of the front face 
$\frf$ of $\bar{X}\x_{\Theta_0}\bar{X}$ such that $\rho_\frf|_{\lb\cap\rb}=\rho_{\pl\frf}$.
Define $\rho_\lb:=\rho/\rho_\frf$ and $\rho_\rb=\rho/\rho_\frf$, they are defining functions 
of $\lb$ and $\rb$, we thus get 
\[[\rho_{\pl\frf}^{4\la}\beta^*((\rho\rho')^{-2\la}R(\la))|_{\lb\cap\rb}]|_{\pl\frf}
=((\rho_\rb\rho_\lb)^{-2\la}\beta^*R(\la))|_{\frf\cap\lb\cap\rb}=\Big((\rho_{\pl\rb}\rho_{\pl\lb})^{-2\la}N(R(\la))\Big)|_{\pl\rb\cap\pl\lb}\] 
where $\pl\rb:=\rb\cap\frf$, $\pl\lb:=\lb\cap\frf$, $\rho_{\pl\rb}=\rho_\rb|_{\frf}$, $\rho_{\pl\lb}=\rho_\lb|_{\frf}$ and $N(R(\la))$ is the normal 
operator of $R(\la)$.
Near a point $p\in M$, we take Darboux coordinates $(u,\omega)\in\rr\x\rr^{2n}$, that is 
$\pl_u=T_0$ and $\Theta_0=du-
(x.dy-y.dx)$ if $\omega=(x,y)\in\rr^n\x\rr^n$, 
they induce coordinates $(\rho,u,\omega;\rho',u',\omega')$
near $(p,p)$ in $\bar{X}\x\bar{X}$. If $T^+\bar{X}=\{d\rho\geq 0\}\subset T\bar{X}$, 
we have the normal fibration $\phi: T^+\bar{X}|_{M}\oplus \rr\pl_{\rho'}\to \bar{X}\x\bar{X}$
\[\phi: \Big(u,\omega;s\pl_\rho+s'\pl_\rho'+tT_0+z.\pl_\omega\Big)
\to \Big(s,u+\demi t,\omega+\demi z; s',u-\demi t, \omega-\demi z\Big)\]
which restricts to $TM\to M\x M$ as a (local near $p$) Riemann-Weyl fibration associated to 
a flat metric, and it induces the fibration $\psi$ since
$T^+\bar{X}|_{M}\oplus \rr\pl_{\rho'}$ is canonically isomorphic to the normal bundle 
$N\pl\Delta(\bar{X}\x\bar{X})$. The fibration $\phi$ lifts to the parabolic
blow-ups around $\pl\Delta$ to give the smooth stucture of $\bar{X}^2_{\Theta}$.
Following Section \ref{normalsection}, the choice of 
homogeneous boundary defining function $\rho_\frf$ of the front face $\frf$
of the blow-up $[N\pl\Delta(\bar{X}\x\bar{X});\pl\Delta,S]$ can be taken to be
\[\rho_\frf:=\Big(\big(s^2+{s'^2}+\demi|z|^2\big)^2+t^2\Big)^{\frac{1}{4}}, \textrm{ with }
|z|^2=d\Theta_0(z',Jz'),\quad z'=z-\Theta_0(z)T_0.\]
Now we observe that 
\[N_p(\Delta_g-\la(n+1-\la))N_p(R(\la))=\delta_{e_p},\]
and the normal operator of $\Delta_g$ is $\Delta_{g_p}$ on 
$(X_p,g_p)$, we deduce from (\ref{resolvent}),(\ref{convolution}) and the group isomorphism $A_p$, 
that the convolution kernel is 
\[N_p(R(\la))\Big(s\pl_\rho+tT_0+\sum_{i=1}^{n}\omega_iX_i'+\omega_{i+n}Y_i'\Big)=c_n\frac{\Gamma(\la)^2}{\Gamma(2\la-n)}r(s,t,\omega)^\la 
{_2F}_1(\la,\la,2\la-n;r(s,t,\omega)),  \] 
\[r(s,t,\omega)=\frac{4s^2}{t^2+(1+s^2+\demi|\omega|^2)^2}=\left(\cosh\left(\frac{d_{g_p}
(s,t,\omega;e_p)}{2}\right)\right)^{-2}\]
with $X_i',Y_i'$ an orthonormal basis of eigenvectors of $J$ for $h_0$ in $\ker\Theta_0$.
Thus we conclude from (\ref{s0}), (\ref{rholbrb}) and 
$\rho_{\frf(Y)}=(t^2+|\omega|^4/4)^{\frac{1}{4}}$ that  
\begin{equation}\label{s_0}
s_0(\la;p,t,\omega)=2c_n(2\la-n-1)\frac{2^{2\la}\Gamma(\la)^2}{\Gamma(2\la-n)}
\left(t^2+\frac{1}{4}|\omega|^4\right)^{-\la}.
\end{equation}
The principal symbol is its Fourier transfom in the fibre $T_p\pl\bar{X}$up Fourier transform of these parabolically homogeneous distributions is studied in great details
by Geller \cite{Ge}, we refer the reader to his work.\\ 

>From Proposition \ref{kernelsla}, Lemma \ref{resvspoisson} and (\ref{noyaup2}) we deduce (see \cite{Gu1,G0} for the proof in the real case)
\begin{lem}\label{scatvspoisson}
For $\la$ such that $P(n-\la),P(\la),S(\la)$ exist, we have the identity
\[P(\la)=P(n+1-\la)S(\la).\]
\end{lem}

Following Ponge \cite[Sec. 5.5]{Po}, one can construct a (``coefficiented'') complex power 
\[L(\la):=c(\la)(1+L)^{\la}\] 
where 
\[L:=\sum_{i=1}^{2n} Z_i^*Z_i, \quad \textrm{span}(Z_1,\dots, Z_{2n})=\ker\Theta_0, \quad c(\la):=\frac{\Gamma(2\la-n-1)}{\Gamma(n+1-2\la)}\] 
then if the metric is even at order $2k$, the operator 
\[\til{S}(\la):=L\Big(-\la+\frac{n+1}{2}\Big)S(\la)L\Big(-\la+\frac{n+1}{2}\Big)\]
is a finite-meromorphic (thus Fredholm) family  in $\cc\setminus (\mc{P}_k\cup (n+1-\mc{P}_k)\cup (n+1-\nn)/2)$  in $\Psi_{\Theta_0}^{0}(\pl\bar{X})$ thus of bounded operators on $L^2(\pl\bar{X})$ by \cite[Prop. 3.1.8]{Po}.  
It is actually holomorphic at $(n+1+\nn)/2$ because of the multiplication by Gamma factors in $\til{S}(\la)$
that kill the first order poles of $S(\la)$ there.
The Heisenberg principal symbol of $p_{2k}=(\la-\la_{2k})S(\la_{2k})-\Pi_{\la_{2k}}$ is given in (\ref{prsymbp2k}) 
and is shown by Ponge \cite[Sec. 3.5]{Po} to be invertible if $\la_{2k}\notin n+1+\nn$, thus
$p_{2k}$ is Fredholm if $\la_{2k}\notin n+1+\nn$. This implies that $\til{S}(\la_{2k})$ is Fredholm on $L^2(\pl\bar{X})$
for these $\la_{2k}$ since $(L+1)^{\la}$ are invertible. Then using Gohberg-Sigal theory \cite{GS} and the fact that $\til{S}(\la)$ is invertible in a small pointed disc centered at $\la_{2k}$, one obtains directly 
from $\til{S}(\la)^{-1}=\til{S}(n+1-\la)$ that $\til{S}(\la)$ is finite meromorphic in a small disc containing 
$n+1-\la_{2k}$ if $n+1-\la_{2k}\notin -\nn_0$ and we can define (as in \cite{G1,PP} for the real case) 
the multiplicity of a pole of finite multiplicity of $\til{S}(\la)$ to be 
\begin{equation}\label{defmults}
\nu(\la_0):=-\tra \textrm{ Res}_{\la=\la_0}(\pl_\la\til{S}(\la)\til{S}^{-1}(\la))=-\tra\textrm{ Res}_{\la=\la_0}
(\pl_\la (c(\la)S(\la))(c(\la)S(\la))^{-1}),
\end{equation}
the second identity being easily checked by cyclicity of the trace (see \cite[Lem. 5.1]{PP}).
Note that $\nu(\la_0)$ is not necessarily finite at $la_0\in-\nn_0$ since $\til{S}(\la)$
is not a priori Fredholm at $n+1-\la_0$ then, and $\til{S}(\la)$ is not a priori finite meromorphic
at $-\nn_0$. Actually it can be proved that it can not be finite meromorphic and that the $p_{2k}$
for $\la_{2k}\in n+1+\nn_0$ have infinite dimensional kernel, this is easy to check for the model $\hh_\cc^{n+1}$. 
This issue and relation with resonances multiplicity at these points will probably be carried out elsewhere. 

Then the method of \cite{G1} based on Gohberg-Sigal theory, Lemma \ref{scatvspoisson}, Lemma \ref{resvspoisson}
and Proposition \ref{kernelsla} can be applied verbatim to obtain 
\begin{prop}\label{multiplicity}
If $g$ is even at order $2k$, if $\la_0\in\{\Re(\la)<(n+1)/2\}\setminus (\mc{P}_k \cup -\nn_0)$ then
we have
\[\nu(\la_0)=m(\la_0)+\indic_{(n+1-\nn)/2}(\la_0)\dim\ker p_{2(n+1-\la_0)} \quad
\textrm{ for }\la_{0}(n+1-\la_0)\notin\sigma_{\rm pp}(\Delta_g)\]
\[ \nu(\la_0)=m(\la_0)-m(n+1-\la_0) \quad \textrm{ for }\la_0\notin \demi(n+1-\nn).\]
\end{prop}
The only non-apparent result we need to apply the proof of \cite[Th. 1.1]{G1} is 
the following unique continuation result dues to Vasy and Wunsch \cite{WV}:
\begin{lem}\label{unique} Let $(X,g)$ be an ACH manifold, and let  $u\in \dot{C}^{\infty}(\bar{X})$ satisfy
 $(\Delta_g-\la(n+1-\la))u=0,$ $\la \in \mr$  in a neighborhood $\Omega$ of $\pl\bar{X},$ 
then $u=0$ in $\Omega.$ In particular, if $\Omega=X,$ $u=0.$
\end{lem}
This is a consequence of Lemma 2.3 of \cite{WV}, after observing that our manifold has bounded geometry  and a product decomposition $[0,\infty)\x \pl\bar{X}$ outside some compact set with metric $g=dr^2+\alpha(r)$ 
for some $1$-parameter of metrics $\alpha(r)$ on $\pl\bar{X}$, sastisfying uniform positivity of second fundamental form
\[-\pl_r\alpha^{-1}(r)>\alpha^{-1}(r),\] 
here $\alpha^{-1}(r)$ is the dual metric to $\alpha(r)$ on $T^*\pl\bar{X}$. We should point out that
the analogue result for conformally compact manifolds was proved by Mazzeo \cite{MA1}.\\

We also remark that the method of \cite{G0} can be used to prove possible essential singularity
for the scattering operator (thus for the resolvent) at $\frac{n-1/2-\nn_0}{2}$:
\begin{prop}\label{esssing}
For each $k\in\nn_0$, there are examples of ACH manifolds for which $n+1-\la_{1+2k}=\frac{n-1/2-k}{2}$ is an essential singularity of $S(\la)$. 
\end{prop}
\textsl{Proof}: Let us consider for simplicity $k=0$ and $\la_1(n+1-\la_1)\notin \sigma_{\textrm{pp}}(\Delta_g)$, 
the other cases are not more complicated (see \cite{G0}), note also that it suffices
to consider $\til{S}(\la)$ instead of $S(\la)$. 
A quick analysis using (\ref{pis}), Subsection \ref{poissonsection} and (\ref{formedelta}) 
shows that the residue of $S(\la)$ at $n+1-\la_1$ is (modulo finite rank smoothing operator)
\[P_1:=\textrm{Res}_{\la=n+1-\la_1}S(\la)=-\frac{n+1-\la_1}{16}(\pl_\rho\log|k(\rho)|)|_{\rho=0}.\] 
Then this operator can be injective on $L^2(\pl\bar{X})$ and any Sobolev space, it suffices actually to 
choose the metric so that this is the Identity (which is easy to do). Then $\til{S}(\la)$
would have a pole of order $1$ whose residue $\til{P}_1$ is injective on $L^2(\pl\bar{X})$.
Suppose that $\til{S}(\la)$ is meromorphic at $n+1-\la_1$, i.e. has a Laurent expansion
at $\la=n+1-\la_1$
\[\til{S}(\la)=\sum_{i=-p}^0S_i(\la-n-1+\la_1)^{i}+O(\la-n-1-\la_1), \quad p\in\nn_0\] 
for some $S_i$. Then since $\til{S}(n+1-\la)\til{S}(\la)=\til{S}(n+1-\la)\til{S}(\la)=\textrm{Id}$ and 
\[\til{S}(\la)=(\la-\la_1)^{-1}\til{P}_1+O(\la-\la_1) \quad \textrm{near }\la_1,\]
we deduce by injectivity of $\til{P}_1$ that $S_i=0$ for $i=-p,\dots,0$. But the term $S_0$
can also be obtained by a contour integral ($C$ is a circle around $n+1-\la_1$ of radius $\eps$)
\[S_0=(2\pi i)^{-1}\int_{C}\frac{\til{S}(\la)}{\la-n-1-\la_1}d\la=\] 
thus this is an operator in $\Psi^{0}_{\Theta_0}(\pl\bar{X})$ with principal symbol
given by 
\[\sigma_{\textrm{pr}}(S_0)=\sigma_{\rm pr}\Big(L\Big(\la_1-\frac{n+1}{2}\Big)\Big)\circ \sigma_{\rm pr}(S(\la))|_{\la=n+1-\la_1}\circ \sigma_{\rm pr}\Big(L\Big(\la_1-\frac{n+1}{2}\Big)\Big)\]
where $\circ$ denotes the composition for Heisenberg principal symbols (see \cite{EMM2,Po}).
But using Theorem \ref{pseudo}
\[\sigma_{\textrm{pr}}(S_0)=\sigma_{\rm pr}\Big(L\Big(-\la_1+\frac{n+1}{2}\Big)\Big)\circ \Big(c_n\frac{2^{n+\frac{3}{2}}\Gamma(\frac{n+1/2}{2})^2}{\Gamma(-\demi)}\mc{F}_{V\to\xi}(||V||_{\rm{He}}^{-2n-1})\Big)\circ \sigma_{\rm pr}\Big(L\Big(-\la_1+\frac{n+1}{2}\Big)\Big),\]
which is not $0$ since the middle one is not $0$, thus a contradiction.
\qed\\

\section{ The Radiation Fields}

In this section we study the scattering theory  developed in the previous sections from a dynamical view point as in the Lax-Phillips theory.  We define the forward and backward radiation fields, show that they give unitary translation representations of the wave group which can be used to define the scattering matrix \eqref{defscat} in terms of the wave equation. The analogue study of  the radiation fields for asymptotically hyperbolic manifolds was carried out in  \cite{SA} .

We start by considering the Cauchy problem for the wave equation
\begin{gather}
\begin{gathered}
\left(D_t^2-\Delta_g-\frac{(n+1)^2}{4}\right) u(t,z)=0  \text{ in } \mr_+ \times X \\
u(0,z)=f_1(z), \;\ D_t u(0,z)=f_2(z), \;\ f_1, f_2 \in C_0^\infty(X).
\end{gathered}\label{we1}
\end{gather}
It is well known that $u\in C^\infty(\mr_+\times X).$ We are interested in understanding the behavior of
$u$ at infinity along some bicharacteristics. 

\begin{theo}\label{radf} Let $\rho$ be a boundary defining function which gives a product decomposition
$\ox=[0,\eps)\times \pl\bar{X}$ for  which \eqref{decompdelta} holds in a collar neighborhood of $\pl\bar{X},$ and denote $\ox \ni z=(\rho,z')\in [0,\eps)\times \pl\bar{X}.$
Let $u(t,z)$ satisfy \eqref{we1}.  Then 
$$v(\rho,s,z'):=\rho^{-n-1}u(s- 2\log \rho, \rho, z)\in C^\infty( [0,\eps)\x \rr \x \pl\bar{X}).$$
\end{theo}
\begin{proof} This is very similar to the proof of Theorem 2.1 of \cite{SA}.    Without loss of generality, we will assume that $f_1=0.$  Let
\begin{gather}
Q=\rho^{-n-2}\left(D_t^2-\Delta_g+\frac{(n+1)^2}{4}\right) \rho^{n+1}. \label{defp}
\end{gather}
Using (\ref{formedelta}) 
and setting  $s=t+2\log \rho,$ we have
\begin{gather}
Q= \p_{\rho} (\p_s+\frac{1}{4} \rho \p_{\rho}) - \rho \Delta_{h(\rho)} +\frac{n+1}{8}A
+\frac{1}{8}A(2\pl_s+\rho\pl_\rho) \label{eqq}
\end{gather}
with $A:=\pl_\rho(\log |k|)$.
Setting $v(\rho,s,z'):=\rho^{-n-1}u(s- 2\log \rho, \rho, z)$, \eqref{we1}  becomes
\begin{gather}
\begin{gathered}
Qv=0. \\
v(\rho, 2\log \rho, z)=0, \;\ \p_s v(\rho,2\log \rho,  z)=\rho^{-n-1} f_2(\rho,z).
\end{gathered}\label{we2}
\end{gather}
>From the standard regularity theory for the wave equation, we know that $v$ is smooth for
$\rho>0.$ The main point of the theorem is to show that $v$ is $C^\infty$ up to $\rho=0.$ 
 There are several ways of proving this, and we will choose a method that will establish
 energy estimates which will be needed in section \ref{supthm}. For that
purpose we make the following change of variables
\begin{gather}
s=4 \log \nu, \;\ \rho=\mu\nu. \label{chv}
\end{gather}
Let 
$$V(\mu,\nu,z)=v(\mu\nu, 4\log \nu, z).$$ In  coordinates  \eqref{chv} equation \eqref{we2} is given by
\begin{gather}
\begin{gathered}
\left( \oq \p_{\mu} \p_{\nu} - \mu\nu \Delta_{h(\mu\nu)} +\frac{n+1}{8}A
+\frac{1}{16}A(\mu\pl_\mu+\nu\pl_\nu)\right)V=0, \\
V(\mu,\mu,z)=0, \;\ (\p_{\mu}V)(\mu,\mu,z)=(\mu)^{-2n-3} f_2(\mu^2, z), \;\
\end{gathered}\label{we3}
\end{gather}
>From Darboux's theorem for contact forms, see for example Theorem 5.5 of \cite{libmar}, we know that for small enough 
$U$ there exists local coordinates $z=(u,x,y),$ $x=(x_1,...,x_n),$ $y=(y_1,...,y_n),$ $u \in \mr,$ near any point 
$z_0\in\pl\bar{X}$ such that
\begin{gather*}
\Theta_0= du +\sum_{j=1}^n \left( y_j dx_j -x_j dy_j\right).
\end{gather*}
Let $X_j$ and $Y_j$ be the vector fields defined by
\begin{equation}\label{basemodel1}
X_j=\frac{1}{\sqrt{2}}(\pl_{x_j}-y_j\pl_u), \quad Y_j=\frac{1}{\sqrt{2}}(\pl_{y_j}+x_j\pl_u), 
\end{equation} 
and $Z_j=X_j$ (resp. $Z_j=Y_{j-n}$) for $j\leq n$ (resp. $j>n$) which forms a basis
of $\ker\Theta_0$. We denote by $T_0=T_0(\rho)$ the Reeb vector field of $(\Theta_0,k(\rho))$.
Note that $\Delta_{h(\mu\nu)}$ is a differential operator on $\pl\bar{X}$
with derivatives in $Z_j,(\mu\nu)^2\pl_u$, since $\rho^2\Delta_{h(\rho)}$ is in $\textrm{Diff}_{\Theta_0}(\bar{X})$.

We want to  differentiate equation \eqref{we3} with  respect to the vector fields 
$X_j$ and $Y_j$ and $\p_u.$  We first recall that
\begin{equation}\label{comm}
[X_j,X_k]=[Y_j,Y_k]=0, \quad [Z_j,\pl_u]=0, \quad [X_j,Y_k]=\delta_{jk}\pl_u.\;\
\end{equation}
Therefore we find that
\begin{gather*}
[X_j^2, Y_j]=2X_j\pl_u, \;\
[Y_j^2,X_j]=-2Y_j\pl_u, \;\ [\pl_u, Z_j^2]=[\p_u^2, Z_j]=0.
\end{gather*}
Differentiating \eqref{we3} and letting $U_1$ be the vector
whose transpose is 
$$U_1^T=(V,X_1V,..., X_n V, Y_1 V, ..., Y_n V, \pl_uV), $$
we deduce that there exists a $(2n+2)\times (2n+2)$ matrix of first order differential 
operator $B$ and functions $C,D,E$ of $(\mu,\nu,z)$ with smooth coefficients such that
\begin{gather}
\begin{gathered}
\Big((\frac{1}{4}\pl_\mu\pl_\nu+\mu\nu\Delta_h)\textrm{Id}  + \mu\nu B(Z_1,\dots,Z_{2n},\mu\nu\pl_u) 
+C\nu\pl_\nu +D\mu\pl_\mu +E\Big)U_p=0, 
\end{gathered}\label{u1}
\end{gather}
and $p=1$.
In general we find that if $U_p$ are the vectors (with size $N=N(p)$) consisting of
all derivatives of $V$ of order $p>1$ with respect to the vector fields
$X_j,$ $Y_j$ and $\frac{\p}{\p u},$ then we get (\ref{u1}) with $B,C,D,E$ 
some $N\x N$ matrices with smooth coefficients.

Notice that the operators in \eqref{u1} are smooth up to $\{\mu=0\}\cup \{\nu=0\}$ and therefore
can be extended smoothly to the neighborhood $\{|\mu|< \sqrt{\mu_0}\}\cup
\{|\nu|< \sqrt{\nu_0}\}.$
 We then prove

\begin{lem}\label{enerest} For $T>0,$ let 
$\Omega=[0,T]\times [0,T]\times \pl\bar{X}.$ 
Let $V(\mu,\nu,y)$ be a $N\times 1$ vector, 
smooth in $\mu>0,$ $\nu>0,$ satisfying the system (\ref{u1})
with 
\begin{gather}
\begin{gathered}
\Big((\frac{1}{4}\pl_\mu\pl_\nu+\mu\nu\Delta_h)\textrm{Id}  + \mu\nu B(Z_1,\dots,Z_{2n},\mu\nu\pl_u) 
+C\nu\pl_\nu +D\mu\pl_\mu +E\Big)V=0, \;\ \text{ and } \\
V(\mu,\mu,y)=f_1(\mu,y), \;\  (\p_\mu V) (\mu,\mu,y)=f_2(\mu,y), \;\ \mu >0,
\end{gathered}\label{eqmcq}
\end{gather}
where $B,C,D,E$ are $N\times N$ matrices with smooth coefficients in $\{|\mu|<T,|\nu|<T,z \in \pl\bar{X}\}$  
with $B$ having for coefficients some first order differential operators in $(Z_j)_j$, $\mu\nu \p_u$. 
If the data  $f_1$ and $f_2$ are such that
\begin{gather}
\begin{gathered}
\int_{0}^T\int_{\pl\bar{X}}\left( \mu | f_1|^2 +
  {\mu}^3|\nabla^hf_1|^2_h 
\right) \; \textrm{dvol}_{k}(z)d\mu 
<\infty, \;\
\int_{0}^T\int_{\pl\bar{X}}
 \mu |f_2|^2 \;\textrm{dvol}_{k}(z)d\mu <\infty,
\end{gathered}\label{condD}
\end{gather}
then for small $T,$ 
\begin{gather}
\begin{gathered}
\int_{ \Omega}\left[ |V|^2+ \mu\nu(\mu+\nu )|\nabla^h V|_h^2+\mu|\pl_\mu V|^2+\nu|\pl_\nu V|^2 \right] \; 
\textrm{dvol}_{k}(z)
d\mu d\nu \leq \\
C\int_{0}^{T}
 \int_{\pl\bar{X}}\left( \mu\left|f_1\right|^2 + \mu\left|f_2\right|^2 
+  {\mu}^3 |\nabla^hf_1|^2_h \right)(\mu,z)\; \textrm{dvol}_{k}(z)
d\mu.
\end{gathered}\label{est3T}
\end{gather}
\end{lem}
\begin{proof} The proof is essentially the same as Lemma 4.1 of \cite{SA}, we just emphasize the differences and we refer the reader
to that paper for details. We denote $Z=(Z_1,\dots,Z_{2n})$ the set of vector fields and we begin by multiplying the system \eqref{eqmcq} by
 $\mu\frac{\p V}{\p \mu}-\nu\frac{\p V}{\p \nu}$, this gives
\begin{gather}
\begin{gathered}
\ha |k|^{-\demi}\pl_\nu\left[ \left( 
\oq \mu\left|\pl_\mu V\right|^2+ {\nu}^2\mu |\nabla^{h} V|^2_h \right)|k|^{\demi}\right]
- \ha |k|^{-\demi}\pl_\mu\left[\left( 
\oq \nu\left|\pl_\nu V\right|^2+ \nu{\mu}^2 |\nabla^{h} V|^2_h\right)|k|^{\demi}
\right] + \\ \sum_k\textrm{div}_k\left(\mu\nu(\mu\p_\mu V- \nu\p_\nu V)\nabla^hV_k
\right)+Q\left(V,\mu\pl_\mu V,\nu\pl_\nu V, \mu\nu ZV, \nu^2\mu^2 \p_u V\right)=0,
\end{gathered}\label{enest1}
\end{gather}
where $Q$ is a quadratic form with smooth coefficients and $\textrm{div}_k$ is the divergence 
with respect to metric $k(\rho)$ on $\pl\bar{X}$, which is easily seen to be also 
the divergence for $h(\rho)$. 

If $(Z_i)_i$ is a local orthonormal basis of $\ker \Theta_0$ for $h(\rho)$ 
(thus for $k(\rho)$) and if $T_0=T_0(\rho)$ is again the Reeb field of $(\Theta_0,k(\rho))$, then
\begin{equation}\label{nablah2}
|\nabla^h V|^2_{h}=\sum_{i}|Z_iV|^2+ \rho^2 a^2(\rho)|T_0V|^2, \quad a(\rho)=1+O(\rho)
\end{equation}
is a smooth function up to $\rho=0$.
Let $\Omega_\del,$ $\Omega_\del^+$  and $\Omega_{a,b},$  be the domains defined by
\begin{gather*}
\Omega_{\del}=[\del,T]\times [\del,T]\times \pl\bar{X}, \;\
\Omega_{\del}^+=\{(\mu,\nu,z)\in \Omega_\del; \;\ \nu \geq \mu\}, \;\ \text{ and }\\
\Omega_{a,b}=\{(\mu,\nu,z)\in \Omega, \;\ a\leq \mu \leq \nu \leq b\}.
\end{gather*}

Integrating in 
$\Omega_{\alpha,T},$ using the compactness of $\p \bar{X},$  the divergence 
theorem, and that the first term in brackets in \eqref{enest1} is positive, we find that for $T$ small, 
there is a constant $C$ such that
\begin{gather}
\begin{gathered}
 \int_{\alpha}^{T} \int_{\p \bar{X}}
\left[\left( 
\nu\left|\pl_\nu V\right|^2+ \nu\mu^2 |\nabla^{h} V|^2_h\right)
\right](\alpha,\nu,z)\; \textrm{dvol}_{k}(z) d\nu 
+\int_{\Omega_{\alpha,T}}|Q|\;\ \textrm{dvol}_{k}(z)d\nu d\mu \\ 
\leq  
C \int_{\alpha}^{T} \int_{\p \bar{X}} \left( 
\mu|\pl_\mu V|^2+\mu|\pl_\nu V|^2
+ {\mu}^3|\nabla^{h} V|^2_h \right)
(\mu,\mu,y)\;\ \textrm{dvol}_{k}(z) d\mu .
\end{gathered}\label{ineq1}
\end{gather}

Doing the same in the region $\Omega_{\delta, \beta}$ gives

\begin{gather}
\begin{gathered}
 \int_{\delta}^{\beta} \int_{\p \bar{X}}
\left[\left( 
\mu\left|\pl_\mu V\right|^2+ \mu\nu^2 |\nabla^{h} V|^2_h\right)
\right](\mu,\beta,z)\; \textrm{dvol}_{k}(z) d\mu 
+\int_{\Omega_{\delta,\beta}}|Q| \;\ \textrm{dvol}_{k}(z)d\nu d\mu \\ 
\leq  
C \int_{\delta}^{\beta} \int_{\p \bar{X}} \left( 
\mu|\pl_\mu V|^2+\mu|\pl_\nu V|^2
+ {\mu}^3|\nabla^{h} V|^2_h \right)
(\mu,\mu,y)\;\ \textrm{dvol}_{k}(z) d\mu .
\end{gathered}\label{ineq2}
\end{gather}

Next we integrate \eqref{ineq1} in $\alpha\in [\delta,T]$ and
\eqref{ineq2} in $\beta\in [\delta,T]$ and add the results to get
\begin{gather}
\begin{gathered}
 \int_{\Omega_\del^+}
\left(  \mu\nu(\mu+\nu)|\nabla^hV|^2_h
+\mu\left|\pl_\mu V\right|^2 
+\nu\left|\pl_\nu V\right|^2 \right)\; \textrm{dvol}_{k}(z) d\mu d\nu + \\
\int_{\del}^{T} \int_{\Omega_{\alpha,T}} |Q| \;\textrm{dvol}_{k}(z)d\mu d\nu d\alpha+ 
\int_{\del}^{T} \int_{\Omega_{\delta,\beta}}|Q| \; \textrm{dvol}_{k}(z) d\mu d\nu d\beta \;\ \leq \\
C(T-\delta)\int_{\del}^{T} \int_{\pl\bar{X}} 
\left( \mu\left|\pl_\mu V\right|^2+ {\mu}^3 |\nabla^{h} V|^2_h \right)
(\mu,\mu,z)\;\ \textrm{dvol}_{k}(z) d\mu.
\end{gathered}\label{firsti} 
\end{gather} 
All terms of $Q,$ except for those in $|V|^2$,  are trivially
bounded by the terms in the first integral, and thus they can be absorbed in first integral
by choosing $T$ small enough. 
To bound the integral of $|V|^2$ in $\Omega_\del^+$ it suffices to copy word by word 
the method in Lemma 4.1 of \cite{SA} and choose $T$ small. This gives 
\begin{gather}
\begin{gathered}
\int_{\Omega_\del^+}|V|^2\; \textrm{dvol}_{k}(z) d\mu d\nu+
\int_{\del}^{T} \int_{\Omega_{\alpha,T}} |V|^2 \;\textrm{dvol}_{k}(z)d\mu d\nu d\alpha+ 
\int_{\del}^{T} \int_{\Omega_{\delta,\beta}}|V| \; \textrm{dvol}_{k}(z) d\mu d\nu d\beta \;\ \leq\\
CT\int_{\Omega_\del^+} \mu\left|\pl_\mu V\right|^2(\mu,\nu,z)\;\ \textrm{dvol}_{k}(z) d\mu d\nu+
C\int_{\delta}^T\int_{\pl\bar{X}}\mu|V|^2(\mu,\mu,z)\;\ \textrm{dvol}_{k}(z) d\mu
\end{gathered} \label{secondi}
\end{gather} 
Now taking $T$ small and $\del\rightarrow 0$ in \eqref{firsti}, \eqref{secondi}, this gives \eqref{est3T} in the region above the diagonal after
using the initial condition. Observing  that the operator in \eqref{eqmcq}
and the estimate \eqref{est3T} remain the same after switching $\mu$ and 
$\nu,$ this estimate also holds in the region below the diagonal. 
This ends the proof of the Lemma.  
\end{proof}

Now we apply the Lemma \ref{enerest}
to prove Theorem \ref{radf}.   The goal is to
prove that for $f_2$ smooth and compactly supported, the solution to
 \eqref{we1}
is smooth up to $\{\rho=0\}.$  We know by  the finite speed  of propagation
that $v$ is
supported in $s\geq C$ for some $C\in\rr$.  The change of variables \eqref{chv}
is smooth in this region and we work in coordinates $(\mu,\nu).$  We will
show that $V(\mu,\nu,z)$ is smooth up to $\{\mu=0\}\cup \{\nu=0\}.$

  We first apply Lemma \ref{enerest} to the special case of
equation \eqref{we3}, noticing that in this case
$V$ is a  function instead of a 
vector. The data on $\{\mu=\nu\}$ is given by \eqref{we3}, so
the last integral on the right hand side of \eqref{est3T} is equal to
\begin{gather*}
\int_{0}^{T} 
\int_{\pl\bar{X}}  \mu|{\mu}^{-2n-3}f_2(\mu^2,y)|^2\;\textrm{dvol}_{k}(z) d\mu \leq \int_{\rho<T} |f_2|^2 \textrm{dvol}_g
\leq ||f_2||_{L^2(X)}^2.
\end{gather*}
Thus the last integral on the right hand side of \eqref{est3T} is bounded
by the square of the norm of $f_2$ in $L^2(X,\dvol_g)$ and
\begin{gather}
\begin{gathered}
\int_{ \Omega}\left[ |V|^2+ \mu\nu(\mu+\nu )|\nabla^h V|_h^2+\mu|\pl_\mu V|^2+\nu|\pl_\nu V|^2 \right] \; 
\textrm{dvol}_{k}(z)
d\mu d\nu \leq 
 C ||f_2||_{L^2(X)}^2.
\end{gathered}\label{l2est}
\end{gather}

Now we want to obtain such energy estimates for the derivatives of $V$. We 
begin by analyzing the derivatives of $V$ with  respect to the vector
fields $X_j,$  $Y_j$ and $\frac{\p}{\p u}$
and we differentiate equation \eqref{we3} with respect to them. 
We  get a system of equations given by \eqref{u1}.
So we conclude in particular that
\begin{gather}
\begin{gathered}
\int_{ \Omega} |U_1|^2+\mu|\pl_\mu U_1|^2+\nu|\pl_\nu U_1|^2+\mu\nu(\mu+\nu)|\nabla^hU_1|_h^2\; \textrm{dvol}_{k}(z)
d\mu d\nu \leq \\
C\left(||f_2||_{L^2(X)}^2+|||\nabla^h f_2|_h||^2_{L^2(X)}\right).
\end{gathered}\label{l2est1}
\end{gather}
Using the commutators \eqref{comm}, the fact that $(Z_j)_j,\pl_u$ form a local basis 
of $T\pl\bar{X}$ near $z_0$,  and that the initial data is smooth
with compact support, \eqref{l2est} and  \eqref{l2est1} guarantee that
\begin{gather*}
V \text{ and } \pl^\alpha_{z} V  \in L^2(\Omega), \; |\alpha|=1.
\end{gather*}
Using this argument repeatedly, we conclude that $V$ is smooth in the variables $z\in\pl\bar{X}$ with values 
in $L^2(\Omega)$ with respect to $(\mu,\nu),$ that is
\begin{gather}
\pl_z^\alpha V\in L^2(\Omega), \;\  \forall \;\ \alpha \in
\mn^{n}, \label{eqq'w}
\end{gather}
and wealso have 
\begin{gather}
\sum_{|\alpha|\leq k}\int_\Omega \left[|\pl_z^\alpha V|^2+\mu|\pl_\mu \pl_z^\alpha V|^2+
\nu|\pl_\nu \pl_z^\alpha V|^2\right] \textrm{dvol}_{k}(z)
d\mu d\nu<\infty, \;\  \forall \;\ k \in
\mn^{n}. \label{eqq''w}
\end{gather}

As in \cite[Th. 2.1]{SA}, we can use the equation that $W:=|k|^{\frac{1}{4}}V$ satisfies to obtain information about 
the derivatives of $V$ with respect to $\mu$ and $\nu$. The proof
is not different of that of Theorem 2.1 in \cite{SA}, we let the reader refer to it
after noticing that Eq. (4.23) there becomes in our case
\begin{gather}
\left(\frac{1}{4}\pl_\nu\pl_\mu +\mu\nu\Delta_{h(\mu\nu)}+ \mu\nu F(\mu\nu)\right)W=0, \label{redeq}
\end{gather}
where $W:=|k|^{\frac{1}{4}}V$ and $F=[ \Delta_{h(\mu\nu)},|k|^\oq]$ is a smooth differential
operator of first order which, in local coordinates, has only derivatives on the vector fields $Z_j,\pl_u$ as above. 
This concludes the proof of Theorem~\ref{radf}.
\end{proof}

With $v$ defined in Theorem \ref{radf}, the map
\begin{gather}
\begin{gathered}
\mcr_+: C_0^\infty(X) \times C_0^\infty(X) \longrightarrow C^\infty(\mr \times \pl\bar{X}) \\
(f_1,f_2) \longmapsto \frac{\p}{\p s} v(0,s,z')
\end{gathered}\label{defrad+}
\end{gather}
will be called the {\it forward} radiation field.  We remark that this definition depends on the choice of $\rho,$ which in turn depends on the choice of the conformal representative in $[\Theta_0]$. The backward radiation field can be defined in a similar
fashion by considering the behavior of the solutions to \eqref{we1} as $t\rightarrow -\infty.$
Let  $u$ be the solution to \eqref{we1}, one can show that
\begin{gather*}
v_-(\rho,s,z')=\rho^{-n-1}u(s+ 2\log \rho, \rho, z)\in C^\infty([0,\eps)\x\rr \times \pl\bar{X}),
\end{gather*}
and define
\begin{gather}
\begin{gathered}
\mcr_-: C_0^\infty(X) \times C_0^\infty(X) \longrightarrow C^\infty(\mr \times \pl\bar{X}) \\
(f_1,f_2) \longmapsto \frac{\p}{\p s} v_-(0,s,z')
\end{gathered}\label{defrad-}
\end{gather}

\section{The Radiation Fields And The Scattering Matrix}\label{Scat}

This section can be compared to section 5 of \cite{SA}: we describe relation between the radiation fields,
the Poisson operator and the scattering operator. As mentioned in subsection \ref{theres}, 
the spectrum of the Laplacian $\sigma(\Delta)$ was studied by Epstein and
 Melrose and Mendoza  \cite{EMM}, and  consists of a finite pure
point spectrum $\sigma_{\pp}(\Delta_g),$ which is the set of $L^2(X)$ 
eigenvalues, and an absolutely continuous spectrum
$\sigma_{\ac}(\Delta_g)$ satisfying
\begin{gather}
\sigma_{\ac}(\Delta_g)= \left [ \frac{(n+1)^2}{4},\infty\right) \;\
\text{ and } \;\ \sigma_{\pp}(\Delta_g)\subset  \left( 0, \frac{(n+1)^2}{4}\right). \label{spect}
\end{gather}

This gives a decomposition 
\begin{gather*}
L^2(X)= L^2_{\pp}(X) \oplus L^2_{\ac}(X),
\end{gather*}
where $L^2_{\pp}(X)$ is the finite dimensional 
space spanned by the eigenfunctions and
$L^2_{\ac}(X)$ is the space of absolute continuity, which is the orthogonal
complement of $L^2_{\pp}(X).$

 The choice of the spectral parameter which adapts well to the wave equation is $\frac{(n+1)^2}{4}+\sigma^2,$ which corresponds to $\la=\frac{n+1}{2}+i\sigma.$   If
$\Im \sigma \not=0,$ then $\frac{(n+1)^2}{4}+\sigma^2 \not \in [\frac{(n+1)^2}{4}, \infty),$
while if $\Im \sigma <-\frac{n+1}{2}$ then
$\frac{(n+1)^2}{4}+\sigma^2 \not \in [0, \infty).$ 
The eigenvalues of $\Delta_g$ are  of finite multiplicity and
are described by points on
the line  $\Re \sigma=0$ and $-\frac{n+1}{2}<\Im \sigma <0.$ The spectral theorem
gives that the resolvent
\begin{gather}
R\Big(\frac{n+1}{2}+i\sigma\Big)=\left(\Delta_g-\frac{(n+1)^2}{4}-\sigma^2\right)^{-1}
: L^2(X)\longrightarrow L^2(X), \;\  \text{ provided } \;\
\Im \sigma <-\frac{n+1}{2}, \label{resb}
\end{gather}
 and it was shown in
\cite{EMM} that it can be meromorphically continued to 
$\mcc \setminus\frac{i}{2}(\ha+\mn_0)$ as an operator acting on appropriate spaces.

Let
\begin{gather*}
H_E(X)=\{(f_1,f_2): \;\ f_1, \; f_2\in L^2(X), \;\ \text{ and } \;
 d f_1 \in L^2(X) \},
\end{gather*}
and for $w_0,  w_1 \in C_0^\infty(X),$ we define
the energy of $w=(w_0,w_1)$ by
\begin{gather}
||w||_E^2=\ha \int_{X} \left( |d w_0|_g^2 -\frac{(n+1)^2}{4}|w_0|^2
+|w_1|^2 \right)\;\ d\vol_g, \label{nore}
\end{gather}
where $|d w_0|_g$ denotes the length of the co-vector with respect to the
metric induced by $g$ on $T^*X.$ 
Note that $||w||_E^2$ is only positive
when $w_0\in L^2_{\ac}(X),$  and only then it defines a norm. 
Let
\begin{gather*}
\Pi_{\ac} : L^2(X) \longrightarrow L^2_{\ac}(X)
\end{gather*}
be the corresponding projector and
\begin{gather*}
E_{\ac}(X)= \Pi_{\ac}\left(H_E(X)\right)=\text{ The range of the projector } 
\Pi_{\ac}  \text{ acting on } \;\ H_E(X),
\end{gather*}
$E_{\ac}(X)$ is a Hilbert space equipped with the norm \eqref{nore}.

Integration by parts shows that if $u(t,z)$ satisfies 
\eqref{we2}, then
$$||\left(u(t,\bullet),D_tu(t,\bullet)\right)||_E=||\left(u(0,\bullet),D_tu(0,\bullet)\right)||_E.$$
The map  $W(t)$ defined by
\begin{gather}
\begin{gathered}
W(t):  C_0^\infty(X)\times C_0^\infty(X)
 \longrightarrow
C_0^\infty(X)\times C_0^\infty(X) \\
W(t)\left(f_1,f_2\right)=\left(u(t,z), D_t u(t,z)\right), \;\ 
 t \in \mr 
\end{gathered}\label{group}
\end{gather}
induces a strongly continuous group of unitary operators that commutes with $\Pi_{\ac}$ 
$$W(t):E_{\ac}(X) \longrightarrow E_{\ac}(X), \;\  t\in \mr.$$

By changing $t\leftrightarrow t-\tau,$ one has that $\mcr_{\pm}$ 
satisfy
\begin{gather}
\mcr_{\pm} \circ \left(W(\tau)f\right)(s,y)=\mcr_{\pm} f(s+\tau,y), 
\;\ \tau \in \mr.
 \label{transl}
\end{gather}

So Theorem \ref{radf}
shows that $\mcr_{\pm}$ are ``twisted'' translation representations of
the group $W(t)$ in the sense of Lax and Phillips.  That is, if one sets 
$\widetilde{\mcr_{\pm}}(f)(s,y)=\mcr_{\pm} f(- s,y),$ then
\begin{gather}
\widetilde{\mcr_{\pm}}( W(\tau))= T_\tau \widetilde{{\mcr_{\pm}}}, \label{transl1}
\end{gather}
where $T_\tau$ denotes right translation by $\tau$ in the $s$ variable. 
Moreover  we will prove
\begin{theo}\label{inv} The maps $\mcr_{\pm}$ induce isometric isomorphisms
\begin{gather*}
\mcr_{\pm} : E_{\operatorname{ac}}(X)\longrightarrow L^2(\mr \times \pl\bar{X}),
\end{gather*} 
where $L^2(\mr \times \pl\bar{X})$ is defined with respect to $h_0$ fixed in 
\eqref{h0} by the choice of boundary defining function $\rho$.
\end{theo}
The proof of Theorem \ref{inv} will be divided into two lemmas. The first one is
\begin{lem}\label{radfl2}  Let $f=(f_1,f_2)\in C_0^\infty(X) \cap E_{\ac}(X).$ Then 
$\mcr_+ f(s,y)\in L^2(\mr\times \pl\bar{X})$ and
\begin{gather*}
||\mcr_+ f||_{L^2(\mr\times \pl\bar{X})}\leq 2 ||f||_E.
\end{gather*}
Moreover, the maps $\mcr_{\pm}$ extend from 
$\left(C_0^\infty(X)\times C_0^\infty(X)\right) \cap E_{\ac}(X)$ by continuity to maps
\begin{gather*}
\mcr_{\pm}: E_{\ac}(X) \longrightarrow L^2(\mr\times \pl\bar{X}).
\end{gather*}
\end{lem}
The proof of this is a straightforward modification of the proof of Proposition 6.4 of \cite{SA} together with
an application of Corollary 6.3 of \cite{SA}.

 Next we have to show that operators $\mcr_{\pm}$ are unitary and onto.  To do that we work on the Fourier transform side and we resort to two observations. The first just follows from taking the partial Fourier transform in the variable $s.$
 \begin{lem}\label{poisson} Let $\mcr_{\pm}$ be the radiation fields defined in \eqref{defrad+}
and \eqref{defrad-}. Then
\begin{gather}
  \begin{gathered}
 \mcf(\mcr_+(f_1,f_2))(\sigma,z)=-\frac{1}{4}\trans{P}\Big(\frac{n+1}{2}+i\sigma\Big)(\sigma f_1-f_2) \\
 \mcf(\mcr_-(f_1,f_2))(\sigma,z)=\frac{1}{4}\trans{P}\Big(\frac{n+1}{2}-i\sigma\Big) (\sigma f_1-f_2).
 \end{gathered} \label{ftradf}
  \end{gather}
where $\mcf$ denotes the partial Fourier transform in the variable $s,$ and $\trans{P}$ is the transpose of the
Poisson operator defined in Proposition \ref{operpoisson}.
\end{lem}
\begin{proof} To see that one needs to observe that
\begin{gather*}
\mcf(\mcr_+(f_1,f_2))(\sigma,z)=-\left.\left(\rho^{-n-1-2i\sigma}i\sigma R(\frac{n+1}{2}+i\sigma)(\sigma f_1-f_2)\right)\right|_{\rho=0} \text{ and } \\
\mcf(\mcr_-(f_1,f_2))(\sigma,z')=\left.\left(\rho^{-n-1+2i\sigma}i\sigma R(\frac{n+1}{2}-i\sigma)(\sigma f_1-f_2)\right)\right|_{\rho=0},
\end{gather*}
This  proof follows that of equation (6.4) of \cite{SA}.  But these operators are the transpose of the forward and backward Poisson operators defined in Proposition
 \ref{poissonprob}.  This proves the lemma.
 \end{proof}

Finally we recall Lemma \ref{resvspoisson} with this spectral parameter
\begin{gather}
\begin{gathered}
R\Big(\frac{n+1}{2}+i\sigma; m,m'\Big)-R\Big(\frac{n+1}{2}-i\sigma ;m,m'\Big)=\\
-4i\sigma
\int_{\pl\bar{X}} P\Big(\frac{n+1}{2}+i\sigma; m,z\Big)P\Big(\frac{n+1}{2}-i\sigma; m',z\Big) \; d\vol_{h_0}(z)
\end{gathered}\label{guiform}
\end{gather}
This and Stone's theorem show that the transpose of Poisson operator
\begin{gather*}
\trans{P} \Big(\frac{n+1}{2}\pm i\sigma\Big): C_0^\infty(X) \longrightarrow C^\infty(\rr^{\pm}\x\pl\bar{X}), \\
\phi \longmapsto \trans{P}\Big(\frac{n+1}{2}\pm i\sigma\Big)\phi(\sigma,z):=c\int_X P 
\Big(\frac{n+1}{2}\pm i\sigma; m,z\Big)\phi(m) \; d\vol_g(m),\;\
\sigma>0
\end{gather*}
extends by continuity to 
\begin{gather}
\begin{gathered}
\trans{P}\Big(\frac{n+1}{2}+i\sigma\Big) : L^2_{\ac}(X) \longrightarrow L^2(\mr^\pm; L^2(\pl\bar{X}), \sigma^2 d\sigma)  \\
\text{ satisfying} \quad
\trans{P}\Big(\frac{n+1}{2} \pm i\sigma\Big)\Delta_g=\Big(\frac{(n+1)^2}{4}+\sigma^2\Big) 
\trans{P}\Big(\frac{n+1}{2} \pm i\sigma\Big)
\end{gathered}\label{specrep}
\end{gather}
and there is $c\in\rr$ such that $c\trans{P}$ gives a surjective isometry between these spaces. 
Lemma \ref{radfl2} is then a straightforward application of these results. See the proof Theorem 5.1 in \cite{SA}.

We deduce from Theorem \ref{inv} that the  {\it dynamical scattering operator}
\begin{gather}
\begin{gathered}
\mcs: L^2(\mr \times \pl\bar{X}) \longrightarrow L^2(\mr\times \pl\bar{X}) \\
\mcs=\mcr_+\circ \mcr_{-}^{-1}
\end{gathered}\label{sco}
\end{gather}
is unitary  in $L^2(\pl\bar{X} \times \mr),$ and in view of \eqref{transl1}, it
commutes with translations. 
  This implies that the Schwartz kernel 
$\mcs(s,y,s',y')$ of $\mcs$ satisfies
\begin{gather*}
\mcs(s,y,s',y')=\mcs\left(s-s',y,y'\right),
\end{gather*}
and thus is a convolution operator.
The {\it scattering matrix} is defined by conjugating $\mcs$ with the partial 
Fourier transform in the $s$ variable
\begin{gather*}
\mca=\mcf \mcs \mcf^{-1}.
\end{gather*}
$\mca$ is a unitary operator in $L^2(\mr \times \pl\bar{X})$ and, since $\mcs$ acts 
as a 
convolution in the variable $s,$ $\mca$ is a multiplication in the 
variable $\la,$ i.e. it satisfies
\begin{gather}
\mca F(\sigma,y)= \int_{\pl\bar{X}} \mca(\sigma;z,z')  F(\sigma,z') \; d\vol_{h_0}(z').
\label{scatmat1}
\end{gather}

We also have
\begin{theo}\label{equiv}  With $\rho$ given by \eqref{modform4} and $\sigma \not=0,$
 the Schwartz kernel of the map $S(\frac{n+1}{2}+i\sigma)$ defined by \eqref{defscat},
 is equal to $-\mca(\sigma;z,z'),$ defined in \eqref{scatmat1}.
\end{theo}
\begin{proof} This follows from \eqref{ftradf} and \eqref{specrep} and the fact that 
\begin{gather*}
S\Big(\frac{n+1}{2}+i\sigma\Big):L^2(\pl\bar{X}) \longrightarrow L^2(\pl\bar{X})
\end{gather*}
is the unitary operator that intertwines $P\left(\frac{n+1}{2}+i\sigma\right)$ and
$P\left(\frac{n+1}{2}-i\sigma\right),$ in the sense that
\begin{gather*}
P\Big(\frac{n+1}{2}+i\sigma\Big)= P\Big(\frac{n+1}{2}-i\sigma\Big)S\Big(\frac{n+1}{2}+i\sigma\Big).
\end{gather*}
\end{proof}

We end this section with a Lemma which will be useful later.
\begin{lem}\label{range1} Let $F\in L^2(\mr \times \pl\bar{X}),$ and let $F^*$ be the function defined
by $F^*(s,z)=F(-s,z),$  then there exists a function $f\in L^2_{\ac}(X)$ such that
$F=\mcr_+(0,f)$ if and only if
\begin{gather*}
F=\mcs F^*.
\end{gather*}
Similarly $F=\mcr_-(0,f)$ if and only if
\begin{gather*}
F^*=\mcs F.
\end{gather*}
\end{lem}
\begin{proof}  Suppose that $F=\mcr_+(0,f).$ Notice that the solution $u$ of the Cauchy problem
 \eqref{we1} with data $(0,f)$ is odd in 
$t.$ Therefore, if $v_+(s,\rho,z)=\rho^{-n-1}u(s-2\log \rho, \rho, z),$ and
$v_-=\rho^{-n-1}u(s+2\log \rho, \rho, z),$ then $v_+(-s,\rho,z)=-v_-(s,\rho,z).$  Hence
$(\p_s v_+)(-s,\rho,z)=\p_s v_-(s,\rho,z).$  Therefore $F(s,z)=\mcr_+(0,f)(s,z)=\mcr_-(0,f)(-s,z).$  This 
implies that $F^*=\mcr_-(0,f)$ and hence $(0,f)=\mcr_-^{-1}F^*.$ So
$F=\mcr_+(0,f)=\mcr_+\mcr_-^{-1}F^*.$  

To prove the converse observe that if $F=\mcr_+(g,f),$ then $F=\mcr_+(g,0)+\mcr_+(0,f)=F_1+F_2.$  If $F=\mcs F^*,$ then $F=\mcs F_1^*+ \mcs F_2^*.$ From the discussion above $\mcs F_2^*=F_2,$ and 
one can prove in a similar way that $F_1=-\mcs F_1^*.$ This shows that $F_1=0$ and by \eqref{specrep} $g=0.$
\end{proof}
\begin{defn}\label{mfmb} With the notation of Lemma \ref{range1} we will denote
\begin{gather}
\mcm^f=\{ F\in L^2(\mr\times \pl\bar{X}): F=\mcs F^*\}, \text{ and }
\mcm^b=\{ F\in L^2(\mr\times \pl\bar{X}): F^*=\mcs F\}. \label{mfmb1}
\end{gather}
\end{defn}

\section{The Support  Theorem} \label{secsupthm}

The goal of this section is to prove
\begin{theo}\label{supthm} Let $f\in L_{\ac}^2(X)$ and suppose that $\mcr_+(0,f)(s,w)=0$ for $s<s_0\ll 0.$ Then $f$ is supported in $\rho>e^{s_0/2}.$
\end{theo}

The proof of this result is divided into two lemmas:
\begin{lem}\label{sup1} Suppose that $f\in L_{\ac}^2(X)$ and that $\mcr_+(0,f)(s,x,y,u)=0$ for $s<s_0\ll 0,$ then there exists $\rho_0>0$ such that $f$ is supported in $\rho>\rho_0.$
\end{lem}
\begin{lem}\label{sup2} Suppose that $f\in C_0^\infty(X)$ and that $\mcr_+(0,f)(s,x,y,u)=0$ for $s<s_0\ll 0,$ then
$f$ is supported in $\rho>e^{s_0/2}.$
\end{lem}

Suppose these two lemmas have been proved.  We now prove Theorem \ref{supthm}
\begin{proof}   We first observe that if  $q\in \mcs(\mr)$ and $\mcf$ is the Fourier transform in $s,$
then it follows from \eqref{specrep}
\begin{gather*}
q(\sigma^2) \mcf \mcr_+(0,f)(\la, z) =\mcf \mcr_+(0, q ( \Delta_g-\frac{(n+1)^2}{4}) f), \;\
f\in L_{\ac}^2(X).
\end{gather*}
If  $\phi\in C_0^\infty(\mr)$ is even then there exists $\psi \in \mcs(\mr)$ such that
$\mcf \phi(\sigma)=\psi(\sigma^2).$ Therefore
\begin{gather}
\phi*\mcr_+(0,f)=\mcr_+(0, \psi(\Delta_g-\frac{(n+1)^2}{4} )f ). \label{phif}
\end{gather}
Since $\p_s^{2k}\phi*\mcr_+(0,f) \in L^2(\mr \times \pl\bar{X}),$ $k=1,2,...,$ then 
\begin{gather*}
(\Delta_g-\frac{(n+1)^2}{4})^k \psi(\Delta_g-\frac{(n+1)^2}{4}) f \in L^2(X).
\end{gather*}
Since $\left(\Delta_g-\frac{(n+1)^2}{4}\right)$ is elliptic in the interior, it follows that 
$\psi(\Delta_g-\frac{(n+1)^2}{4})f \in C^\infty(X).$

Let $f\in L^2_{\ac}(X)$  be such that  $F(s,w)=\mcr_+(0,f)(s,w)=0$ for $s<s_0\ll 0.$ 
Let $\phi\in C_0^\infty(\mr)$ be even and supported in
$|s|<1$ and $\int \phi(s) \; ds=1.$  Let $\phi_\eps(s)=\eps^{-1}\phi(s/\eps)$ and $F_\eps=\phi_\eps*F,$
then $F_\eps$ is supported in $s>s_0-\eps.$ On the other hand, in view of \eqref{phif} there is
 $f_\eps\in L_{\ac}^2(X)$ such that $\mcr_+(0,f_\eps)=F_\eps.$ In view of Lemma \ref{sup1}
 $f_\eps$ is compactly supported.  On the other hand, since $f_\eps$ is smooth, Lemma \ref{sup2}
 guarantees that $f_\eps$ is supported in $\rho \geq e^{\frac{s_0-\eps}{2}}.$  Letting $\eps\rightarrow 0$
 it follows that $f$ is supported in $\rho\geq e^{s_0/2}.$
\end{proof}
Now we prove Lemma \ref{sup2}.
\begin{proof}    Suppose the initial data is supported in $\rho\geq \rho_0.$  By finite speed of propagation the solution to \eqref{we2} is equal to zero in the region
\begin{gather}
\{4\log \rho-2\log \rho_0 \leq s \leq 2\log \rho_0\}.\label{vanish}
\end{gather}
We have shown that the solution to Equation \ref{we2} with compactly supported initial data is smooth up to $\rho=0.$ By assumption we know that $\frac{\p v}{\p s}(0,s,z)=0$ if $s\leq s_0.$  One can then show that this implies that $v(\rho,s,y)$ vanishes to infinite order at $\rho=0$ provided $s\leq s_0.$  Therefore we deduce from \eqref{vanish} that the solution $v$ to \eqref{we2} vanishes to infinite order at 
the corner $\{\rho=0\}\cup\{s=2\log \rho_0\}.$ In particular we conclude that $v$ vanishes to infinite order at
$\rho+s=2\log \rho_0$ and is supported in the wedge  $\{\rho\geq 0\}\cup\{s\geq \log \rho_0\}.$ Now we want to use
a uniqueness theorem that would allow us to show that $v=0$ in a neighborhood of the set
$\{s=\log \rho_0,\rho=0\}\times \pl\bar{X}.$   The result we need is a particular case of a Theorem due to 
Alinhac, Theorem 1.1.1 of \cite{ALI}. We will verify that the hypotheses of this theorem are satisfied in this case.  

The principal symbol of the operator $Q$ in \eqref{we2} is
\begin{gather*}
q=-\tau\sigma -\oq \rho \tau^2-\rho p
\end{gather*}
where $(\tau,\sigma,\xi,\eta,\nu)$ is the dual to $(\rho,s,x,y,u),$ and $p=p(\rho,x,y,u,\xi,\eta,\nu),$ is the symbol of the Laplacian
$\Delta_{h(\rho)}.$   It follows from \eqref{decompdelta} and \eqref{basemodel1} that one can choose local coordinates $(x,y,u)$ near a point on $z_0=(x_0,y_0,u_0) \in \p X$ such that
\begin{gather}
p=\sum h^{ij}(\xi_i-y_i\nu, \eta_j+x_j\nu)+\rho^2 \nu^2 + \rho^2 p_1 +O(\rho^3),\label{locp}
\end{gather}
 where $p_1$ does not depend on $\nu,$ and
$h^{ij}$ is a positive definite matrix. It is important to observe that $p$ is elliptic if $\rho\not=0.$

The Hamilton vector field of $q$ is
\begin{gather*}
H_q= -\tau\p_s-\sigma \p_\rho -\ha\rho\tau\p_s+
\left(p+\rho\p_\rho p +\oq \tau^2\right)\p_\tau -\rho H_p,
\end{gather*}
where $H_p$ is the Hamilton vector field of $p$ with respect to the tangential variables, and it does not have derivatives in $s$ and $\rho,$ $\sigma$ or $\tau.$

Let $\phi=s+\rho-2\log \rho_0,$  and let $\{f,g\}=H_f g$ be the Poisson bracket between $f$ and $g,$ then 
\begin{gather}\label{pphi}
\{q,\phi\}=H_q \phi=-\tau -\sigma -\ha\rho\tau.
\end{gather}
So, for $\rho$ small,
\begin{gather}
\{\{q,\phi\},\phi\}= H_{\{q,\phi\}} \phi= -2-\ha \rho \not= 0 \label{ppphi}
\end{gather}
We also have
\begin{gather*}
H_q^2 \phi=\ha \sigma\tau+\frac{1}{8}\rho\tau^2 -\oq \tau^2 -(1+\ha \rho)(p+\rho\p_\rho p).
\end{gather*}
Then, for $\rho$ small,
\begin{gather}
\; q=H_q\phi=0  \Longrightarrow  H_q^2\phi=-\oq \tau^2-(1+\rho) p -(1+\ha\rho) \rho\p_\rho p\leq 0
\label{neg}
\end{gather}
since $-\rho\pl_\rho p\leq p$, which is a consequence of 
\[-\rho\pl_\rho (\rho^2\nu^2)=-2\rho^2\nu^2\leq 0 , \quad |\rho\pl_\rho p_1|\leq p\]
for small $\rho$.
We want to analyze the set 
$\Lambda=\{q= H_q\phi=H_q^2\phi=0\}\subset T^*(\rr\x\bar{X})$ near $\rho=0.$
Combining the equations defining $\Lambda,$ and the ellipticity of $p$ for $\rho>0,$ one concludes that, for small $\rho,$
\begin{gather}
\Lambda=\{q= H_q\phi=H_q^2\phi=0\}
=\{\rho=\sigma=\tau=p=0\}. \label{Lambda}
\end{gather}
In local coordinates near a point $z_0\in \p X$ where \eqref{locp} is valid
\begin{gather}
\Lambda=\{\rho=\sigma=\tau=0, \; \xi_j=y_j\nu, \; \eta_j=-x_j\nu\}, \label{Lambda1}
\end{gather}
and therefore it is a smooth submanifold. It also follows from \eqref{locp} that
\begin{gather}
dp=0 \;\ \text{ on } \;\ \Lambda. \label{dp=0}
\end{gather}

One can check that the symbol $e$ given by
\begin{gather}
e=-\ha H_q^2 \phi+\oq(H_q\phi)^2 \geq 0  \text{ and is transversally elliptic to } \Lambda.\label{symb}
\end{gather}
Finally, notice that
\begin{gather*}
H_{\{q,\phi\}}=-(1+\ha \rho)\p_{\rho} -\p_s-\ha\tau\p_{\tau},
\end{gather*}
and so
\begin{gather}
H_{\{q,\phi\}} \;\ \text{ is transversal to } \;\ \Lambda. \label{transvto}
\end{gather}

 So the following conditions are satisfied near any point 
 $(0,2\log \rho_0,z_0) \in \{\rho=0, s=2\log\rho_0\}:$    \eqref{ppphi}, \eqref{neg} hold,
 $\Lambda$ is a  smooth submanifold, \eqref{dp=0}, \eqref{symb}
 and  \eqref{transvto} are true. Moreover,  $v$ is supported in the wedge $\{\rho\geq 0\}\cap \{ s \geq 2\log\rho_0\},$ so
 the intersection of $\{\phi=0\}$ with the support of $v$ is compact.
   Then it follows from Theorem 1.1.1 of \cite{ALI} that $v=0$ in a neighborhood  of  the point  
   $(0,2\log \rho_0,x_0,y_0,u_0).$ Using the compactness of $\p \bar{X},$
   we conclude that $v=0$ in a neighborhood of the set $\{\rho=0, s=2\log \rho_0\}.$

Repeating this argument one concludes that there exists $\del>0$ such that $v(\rho,s,z)=0$ if
$\rho<\del$ and $2\log \rho_0\leq s \leq s_0.$   Now we can repeat the argument starting at 
$\{\rho=\del, s=2\log \rho_0\}.$  This is actually easier, because the level surfaces of 
$\phi=s+\rho$ are strictly pseudoconvex away from $\rho=0.$ Therefore H\"ormander's Theorem, see Theorem 28.3.5 of \cite{HO}, can be used to show that $v=0$ in a neighborhood of $\{\rho=\del,s=2\log \rho_0\}.$
This process can be continued to show that $v=0$ in the region 
$\{0\leq \rho\leq \rho_0, 2\log \rho_0\leq s \leq s_0\}.$  Now we translate this back to the $t$ variable, and
using that $u$ is odd in $t,$ we conclude that $u(t,z)=0$ in  $\rho\leq \rho_0$ and
$0 \leq t \leq s_0-2\log \rho_0.$  Now Tataru's theorem in \cite{TAT} shows that
$\p_t u(0,z)=f_2(z)=0$ if $2\log \rho\leq s_0.$  This proves the  Lemma
\end{proof}

Now we prove Lemma \ref{sup1}.
\begin{proof}    To extract information about the behavior of $v$ as $\rho \rightarrow 0$ and $s\rightarrow-\infty$ we need to work with the compactified equation \eqref{we3}.  In fact we will work with 
\eqref{redeq} which does not have first order derivatives in $\mu$ and $\nu.$  This is very similar to the proof of Lemma 7.2 of \cite{SA}.  

By \eqref{phif}, we can assume that the initial data $f \in C^\infty(X).$ Therefore, by standard regularity
for solutions to the wave equation, the solution $W$ to \eqref{redeq} is smooth in the region $\{\mu>0, \nu>0\}.$ 
First we will show that,
as a distribution,  $W$  vanishes to infinite order at $\{\mu=0\}\cup\{\nu=0\}.$  If we knew that $W$ were smooth up to  $\{\mu=0\}\cup\{\nu=0\},$ then  Theorem 1.1.2 of \cite{ALI}
would guarantee that $W=0$ near $\{\mu=\nu=0\}$ and this would imply in particular that the initial data
is supported away from $\rho=0.$  However, we do not know this in principle, and we will use the fact that $\mcr_+(0,f)=0$ in $s<s_0$ to show that this is true in a neighborhood of $\{\mu=\nu=0\}.$

The first step is to show that $W$  which is defined in $(0,T)_\mu\times (0,T)_\nu \times \pl\bar{X}$, has an extension 
\begin{gather*}
\tW\in H^{2k}((-T,T)\times (-T,T); H^{-2k}(\pl\bar{X})), \;\ \text{ for all } k \in \mn,
\end{gather*}
which satisfies \eqref{redeq}.  The proof of this fact is identical to the proof of the analogous case done in \cite{SA},  the only ingredient needed here are the energy estimates from Lemma \ref{enerest}.

Since we cannot apply Alinhac's theorem directly, we need to verify that the regularity of $V$ given by
\eqref{enerest} is enough for the methods of \cite{ALI} to work.   We first prove the a Carleman estimate.  Let $u$ be a smooth function of its arguments that vanishes to infinite order at $\{\mu=0\}\cup\{\nu=0\}.$
 Let
\begin{gather}
Pu=\left(\frac{1}{4}\pl_\nu\pl_\mu +\mu\nu\Delta_{h(\mu\nu)}+ \mu\nu F(\mu\nu)\right)u \label{redeq1}
\end{gather}

Let $u=(\mu+\nu)^\gamma v$ and let  $P_\gamma v=(\mu+\nu)^{-\gamma} P (\mu+\nu)^{\gamma} v.$ Then
\begin{gather*}
P_\gamma v= P v+ \frac{\gamma(\gamma-1)}{4(\mu+\nu)^2} v +
\frac{\gamma}{4(\mu+\nu)}(\p_\mu  + \p_\nu)v
\end{gather*}
Let $\Omega=[-T,T]\times[-T,T]\times \pl\bar{X},$ let  $\lan f, g \ran,$ denote the $L^2(\Omega)$ 
inner product of $f$ and $g,$ and let $||f||=||f||_{L^2(\Omega)}.$ Therefore
\begin{gather*}
||P_\gamma v||^2= ||(P+  \frac{\gamma(\gamma-1)}{4(\mu+\nu)^2} v||^2+
||\frac{\gamma}{4(\mu+\nu)}(\p_\mu  + \p_\nu)v||^2 \\ +
2\lan (P+  \frac{\gamma(\gamma-1)}{4(\mu+\nu)^2} )v,
\frac{\gamma}{4(\mu+\nu)}(\p_\mu v + \p_\nu)v\ran.
\end{gather*}
Let $\Sigma_1=\{\mu=T\},$ $\Sigma_2=\{\nu=T\}.$ Using Hardy's inequality
$||(\mu+\nu)^{-1}(\p_\mu+\p_\nu)v||^2\geq \left(9/4\right)||(\mu+\nu)^{-2}v||^2,$
and integrating by parts we find that for $T$ small there exists a constant $C>0$ such that
\begin{gather}
\begin{gathered}
||P_\gamma v||^2 +C\gamma \int_{\Sigma_1} \frac{\mu\nu}{\mu+\nu}\left |\nabla^h v\right| \; d\sigma 
+C\gamma \int_{\Sigma_2} \frac{\mu\nu}{\mu+\nu}\left |\nabla^h v\right| \; d\sigma \geq  
 \gamma^2 ||(\mu+\nu)^{-1}(\p_\mu +\p_\nu)v||^2 \\ +\gamma ||(\mu+\nu)^{-1}\p_\mu v||^2 + 
 \gamma ||(\mu+\nu)^{-1}\p_\mu v||^2 + \gamma^2 ||(\mu+\nu)^{-2} v||^2 + ||\nabla^h v||^2.
 \end{gathered}\label{carle1}
 \end{gather}
 Using that $v=(\mu+\nu)^{-\gamma }u,$ we get
 \begin{gather}
 \begin{gathered}
 ||(\mu+\nu)^{-\gamma} P u||^2 +C\gamma \int_{\Sigma_1} \mu\nu(\mu+\nu)^{-\gamma-1}
 \left |\nabla^h v\right| \; d\sigma 
+C\gamma \int_{\Sigma_2} \mu\nu(\mu+\nu)^{-\gamma-1}\left |\nabla^h v\right| \; d\sigma \geq  \\
 \gamma^2 ||(\mu+\nu)^{-1}(\p_\mu +\p_\nu)(\mu+\nu)^{-\gamma}u||^2 +
 \gamma ||(\mu+\nu)^{-1}\p_\mu (\mu+\nu)^{-\gamma-1}u||^2 + \\
 \gamma ||(\mu+\nu)^{-1}\p_\mu (\mu+\nu)^{-\gamma-1}u ||^2 + 
 \gamma^2 ||(\mu+\nu)^{-\gamma-2} u||^2 + ||(\mu+\nu)^{-\gamma} \nabla^h u||^2.
 \end{gathered}\label{carle2}
 \end{gather}
 
 Now we want to apply this to $\tW$ which is the solution to \eqref{redeq}. To do that we have to regularize $\tW.$  Let $z_0\in \pl\bar{X}$ and let $U_0$ be a neighborhood of $z_0.$  Let 
 $\psi\in C_0^\infty(U_0)$ with $\psi(z_0)=1.$ Let $\chi\in  C_0^\infty(\mrn),$ $\chi(0)=1$ and
 $\int \chi=1.$  Let $\chi_\del(w)=\del^{-n}\chi(w/\del),$ and let $u_\del=\chi_\del*(\psi \tW).$  Then
 we can apply \eqref{carle2} to $u_\del.$ Now we want to let $\del\rightarrow 0.$  This is done, as usual,
 by using Friedrich's lemma, see for example Theorem 2.4.3 of \cite{HO1}.   To do that we need to know that the right hand side of \eqref{carle2} is finite for some $\gamma>0.$ 
 Since $P\tW=0,$ we can use estimates \eqref{enerest},
 \eqref{ineq1} and \eqref{ineq2} to deduce that the left hand side of \eqref{carle2} is finite for 
 $\gamma=\ha.$  So we deduce that \eqref{carle2} holds with $\gamma=1/2$ and $u=\psi \tW.$ This implies that we can use Friedrich's lemma when $\gamma=3/2,$ and so \eqref{carle2} holds for 
 $u=\psi\tW$ and $\gamma=3/2.$  Then a bootstrapping argument shows that it holds for all $\gamma.$  So we conclude that  \eqref{carle2} holds for $\psi\tW,$  and all $\gamma.$  Now a standard partition of unity argument gives $\tW=0$  in $\Omega$ for $T$ small.
\end{proof}

The following are  consequences of the support theorem which is the key step in the reconstruction
of the manifold from the scattering matrix.
\begin{lem}\label{range2} Let  $\mcr_{\pm}$ be the forward an backward radiation fields defined above.  Let $a>0$ be small.  Then
\begin{gather*}
\mcm^f(a):=\\ \{ \mcr_{+}(0,f); f \in L^2(X), \; f=0 \; \text{ if } \rho<a\}=\{F\in L^2(\mr\times \pl\bar{X}): F=\mcs F^* ,
F=0 \text{ if } s< 2\log a\} \\
\mcm^b(a):= \\ \{ \mcr_{-}(0,f); f \in L^2(X), \; f=0 \; \text{ if } \rho<a\}= 
\{F\in L^2(\mr\times \pl\bar{X}): F^*=\mcs F,
F=0 \text{ if } s>- 2\log a\} \\
\end{gather*}
\end{lem}
\begin{proof}  From Lemma \eqref{range1} $F=\mcr(0,f)$ if and only if
$F=\mcs F^*.$ The support theorem guarantees that $F=0$ for $s<2\log a$ if and only if $f$ is supported
in $\rho\geq a.$
\end{proof}
The following lemma will be very important in the reconstruction of the manifold
\begin{lem}\label{proj} Let $\mcm^f$ and $\mcm^b$ be the spaces defined in \eqref{mfmb1}.  Let
\begin{gather}
\begin{gathered}
\Pi^f(a): \mcm^f \longrightarrow \mcm^f(a)\\
\Pi^b(a): \mcm^b \longrightarrow \mcm^b(a)
\end{gathered}\label{proj1}
\end{gather}
be the orthogonal projections.  Then for any $a\in \mr,$ $\Pi^f(a)$ and $\Pi^b(a)$ are determined by the scattering matrix. 

Moreover, there exists $\del>0$ such that if $a \in [0,\del)$ and if $\chi_a$ is the characteristic function of the set $\{\rho\geq a\},$ then, if   
$(\bar{X},g) $ has no eigenvalues,
\begin{gather}
\Pi^f(a) \mcr_+(0,f)= \mcr_+(0, \chi_a f), \;\ 
\Pi^b(a) \mcr_+(0,f)=\mcr_-(0,\chi_a f).
\end{gather}
If   $(\bar{X},g)$  has eigenvalues there exist  a continuous family of finite rank operators
$T(a)$ such that
\begin{gather}
\Pi^f(a) \mcr_+(0,f)=\mcr_+(0, \chi_a (I+T(a)) f), \;\
\Pi^b(a) \mcr_-(0,f)=\mcr_-(0, \chi_a (I + T(a) )f). \label{proj2}
\end{gather}
\end{lem}
\begin{proof} Since the spaces $\mcm^f,$ $\mcm^f(a),$ $\mcm^b$ and $\mcm^b(a)$ are determined by the scattering matrix, the first part  is immediate. 

Since $\Pi^f(a)\mcr_+(0,f)$ is supported in $s\geq 2\log a,$ the support theorem guarantees that there exists $f_a\in L^2_{\ac}(X)$ supported in $\rho \geq a$ such that
\begin{gather*}
\Pi^f(a)\mcr_+(0,f)=\mcr_+(0,f_a).
\end{gather*}
Let $g \in L^2_{\ac}(X)$ supported in $\rho \geq a,$ then
\begin{gather*}
\lan \Pi^f(a) \mcr_+(0,f), \mcr_+(0,g)\ran = \lan f_a, g\ran =\lan f, g\ran.
\end{gather*}
Then
\begin{gather}
\lan f_a-f,g\ran=0 \;\ g \in L^2_{\ac}(X), \text{ supported in } \rho \geq a. \label{ortho}
\end{gather}

 If $(\bar{X},g)$ has no eigenvalues, then
$L^2_{\ac}(X)=L^2(X).$  In this case we conclude that $f_a=\chi_a f.$

When there are eigenfunctions we can only deduce that in $\rho\geq a,$ $f_a-f=\sum_j c_j(a,f) \phi_j.$ In other words,
\begin{gather*}
f_a=\chi_a( f + \sum_{j=1}^n c_j(a,f) \phi_j).
\end{gather*}
Since $f_a\in L^2_{\ac}(X),$ $\lan f_a, \phi_j\ran =0,$  and so we have
\begin{gather}
\lan \chi_a f, \phi_k\ran + \sum_{j=1}^N c_j(a,f) \lan \chi_a \phi_j, \phi_k\ran =0. \label{system}
\end{gather}

Since $\{\phi_j\}_{j=1}^N$ is orthonormal, the system \eqref{system} can be solved for $a=0.$ Therefore there exists $\del>0$ such that it can be solved for $a \in [0,\del].$   It is clear that the constants $c_j$ depend linearly  on $f.$ Therefore the operator $f \longmapsto \sum_{j=1}^Nc_j(a,f) \phi_j$ is linear and of finite rank.  This is the operator $T(a)f.$
\end{proof}

\section{ The Inverse Problem}

    In this section we will prove Theorem \ref{inverseprob}. The proof is based on the boundary control theory of Belishev, and the key point is the support theorem proved in section \ref{secsupthm}.
     
First note that the result in Theorem \ref{pseudo} about the principal symbol
of $S(\la)$ implies that the pseudo-hermitian structure on $\pl\bar{X}$ are the same for $(X,g_1)$ and $(X,g_2)$
if $S_1(\la)=S_2(\la)$ and $S_1(\la),S_2(\la)$ defined from the same conformal representative of $[\Theta_0]$.
Thus the dilation $M_\rho$ defined in Subsection \ref{ACHM} are the same for both metrics, as well as 
the metrics $h_0=\Theta_0^2+d\Theta_0(.,J.)$. 

We begin with the following lemma:  
\begin{lem}\label{inv1} Let $(X_j,g_1)$ be ACH manifolds satisfying the hypotheses of Theorem \ref{inverseprob}.
For $j=1,2$, there exist diffeomorphisms $\psi_j:[0,\eps)_\rho\x M\to \psi_j([0,\eps)\x M)\subset\bar{X}_j$ with 
$\psi_j(\{\rho=0\})=M$ and 
\begin{gather*}
\psi_j^*g_j=\frac{d\rho^2+h_j(\rho)}{\rho^2}, \quad M_\rho^*(\rho^{-2}h_j(\rho))|_{\rho=0}=h_0
\end{gather*}
for some smooth family of metrics $h_j(\rho)$ on $M$ in $\rho\in(0,\eps)$ and such that $M_\rho^*(\rho^{-2}h_j(\rho))$ are smooth metrics on $M$
depending smoothly on $\rho\in[0,\eps)$. 
\end{lem}
\begin{proof} This is just following Subsection \ref{modelbdfs}.
\end{proof}
The next step is to prove
\begin{prop}\label{h1h2}  For $j=1,2$, let $(X_j,g_j)$ be ACH manifolds satisfying the hypotheses of 
Theorem \ref{inverseprob}, and let $h_j(\rho)$ given by Lemma \ref{inv1}.
Then there exists $\del\in (0,\eps)$ such that $h_1(\rho)=h_2(\rho)$ for $\rho \in  [0,\del).$
\end{prop}
\begin{proof}    Let $F\in \mcm^b$  be such that $\p_s^{2k} F \in L^2(\mr\times \pl\bar{X}),$ for all $k \in \mn.$  Let $a\in \mr,$ $a\ll 0.$  Let $F=\mcr_+(0,f).$ Then, $f$ is $C^\infty,$ and  by Lemma \ref{proj}
\begin{gather*}
\mcr_-^{-1}\Pi^b(a) F=(0, \chi_a f), \text{ if there are no eigenvalues} \\
\mcr_-^{-1}\Pi^b(a) F=(0, \chi_a (I+T(a))f), \text{ if there are  eigenvalues}.
\end{gather*}

Therefore
\begin{gather}
\begin{gathered}
\mcr_+\mcr_-^{-1}\Pi^b(a) F=\mcs\Pi^b(a)F =\mcr_+(0, \chi_a f), \text{ if there are no eigenvalues} \\
\mcr_+\mcr_-^{-1}\Pi^b(a) F=\mcs\Pi^b(a) F=\mcr_+(0, \chi_a (I+T(a))f), \text{ if there are eigenvalues} 
\end{gathered}\label{proj3}
\end{gather}
The left hand side of  each equation of  \eqref{proj3} is determined by the scattering matrix. Therefore so is the right hand side.  The initial data is singular at $\rho=a,$ and this singularity will travel to the boundary as $t\rightarrow \infty.$ We want  to find the singularity of 
$\mcr_+(0,\chi_a f)$ or $\mcr(0, \chi_a(I+T(a))f)$ at $s=2\log a.$ 

This can be done exactly as in the proof of Lemma 8.9 of \cite{SA}, and we find that for $a\in [0,\del/4],$
with $\del$ given by Lemma \ref{proj}
\begin{gather}
\begin{gathered}
\mcs\Pi^b(a)F(s,z) =\mcr_+(0, \chi_a f)(s,z)= \ha a^{-n-1} f(a,z) |k|^\oq(a,z) |k|^{-\oq}(0,z)+ \\
\text{smoother terms if there are no eigenvalues} \\
\mcs\Pi^b(a)F(s,z) =\mcr_+(0, \chi_a(I+T(a)) f)(s,z)= \ha a^{-n-1}(I+T(a)) f(a,z) |k|^\oq(a,z) |k|^{-\oq}(0,z)+
\\ \text{smoother terms if there are eigenvalues} 
\end{gathered}\label{proj4}
\end{gather}
 
 Doing the same with $F$ replaced by $\p_s^2F$ we find that
 \begin{gather}
\begin{gathered}
\mcs\Pi^b(a)\p_s^2 F(s,z) =\mcr_+(0, \chi_a (\Delta_g-\frac{(n+1)^2}{4}) f)(s,z)= \\
\ha a^{-n-1}\left((\Delta_g- \frac{(n+1)^2}{4}) f\right)(a,z) |k|^\oq(a,z) |k|^{-\oq}(0,z)+ \\
\text{smoother terms if there are no eigenvalues} \\
\mcs\Pi^b(a)\p_s^2 F(s,z) =\mcr_+(0, \chi_a(I+T(a)) (\Delta_g-\frac{(n+1)^2}{4})f)(s,z)= \\ \ha a^{-n-1}(I+T(a))\left((\Delta_g-\frac{(n+1)^2}{4}) f\right)(a,z) |k|^\oq(a,z) |k|^{-\oq}(0,z)+
\\ \text{smoother terms if there are eigenvalues} 
\end{gathered}\label{proj5}
\end{gather} 
Now consider the manifolds $(X_1,g_1)$ and $(X_2,g_2)$ satisfying the hypotheses of Theorem \ref{inverseprob}.  Let $(0,f_j)=\mcr_{-,j}^{-1}F.$ Since \ $\mcs_1=\mcs_2,$ then, if there are no eigenvalues,
\begin{gather*}
f_1(a,z) |k_1|^\oq(a,z) |k_1|^{-\oq}(0,z)= f_2(a,z) |k_2|^\oq(a,z) |k_2|^{-\oq}(0,z), \text{ and } \\
\left((\Delta_{g_1}- \frac{(n+1)^2}{4}) f_1\right)(a,z) |k_1|^\oq(a,z) |k_1|^{-\oq}(0,z)\\=\left((\Delta_{g_2}- \frac{(n+1)^2}{4}) f_2\right)(a,z) |k_2|^\oq(a,z) |k_2|^{-\oq}(0,z) 
\end{gather*}

Therefore substituting the first equation in the second, we get
\begin{gather*}
\left((\Delta_{g_1}- \frac{(n+1)^2}{4}) f_1\right)(a,z) |k_1|^\oq(a,z) |k_1|^{-\oq}(0,z)\\=\left((\Delta_{g_2}- \frac{(n+1)^2}{4}) \left( f_1 |k_1|^\oq |k_2|^{-\oq}\right) \right)(a,z) |k_2|^\oq(a,z) |k_2|^{-\oq}(0,z)
\end{gather*}
implies that the operators have the same coefficients,
and hence $|k_1|=|k_2|$ and $\Delta_{g_1}=\Delta_{g_2}.$

When there are eigenvalues, we have the following
\begin{gather*}
(I+T_1(a)) f_1(a,z)  |k_1|^\oq(a,z) |k_1|^{-\oq}(0,z)= 
(I+T_2(a)) f_2(a,z) |k_2|^\oq(a,z) |k_2|^{-\oq}(0,z), \text{ and } \\
(I+T_1(a))\left((\Delta_{g_1}-\frac{(n+1)^2}{4}) f_1\right)(a,z)
 |k_1|^\oq(a,z) |k_1|^{-\oq}(0,z)= \\
(I+T_2(a))\left((\Delta_{g_2}-\frac{(n+1)^2}{4}) f_2\right)(a,z) |k_2|^\oq(a,z) |k_2|^{-\oq}(0,z).
\end{gather*}
Hence, proceeding as in the case of no eigenvalues, we obtain
\begin{gather*}
\left((\Delta_{g_1}- \frac{(n+1)^2}{4}) f_1\right)(a,z) |k_1|^\oq(a,z) |k_1|^{-\oq}(0,z)\\=\left((\Delta_{g_2}- \frac{(n+1)^2}{4}) f_1\right)(a,z) |k_2|^\oq(a,z) |k_2|^{-\oq}(0,z) + Tf(a,z),
\end{gather*}
where $T$ is an operator of finite rank.

The difference between the term on the left hand side of this equation and the first term of the right hand side is a differential operator, while the second term on the right hand side is an operator of finite rank.  Therefore, the differential operators must be equal. So $T=0,$ and we argue as above to
conclude that $|k_1|=|k_2|$ and $\Delta_{g_1}=\Delta_{g_2}.$
\end{proof}

Proposition \ref{h1h2} shows that there exist $\eps>0$ and  a smooth diffeomorphism
\begin{gather*}
\Psi : M \times [0,\eps) \longrightarrow M\times [0,\eps) \\
\Psi^*(g_2)=g_1.
\end{gather*}

To extend the diffeomorphism $\Psi$ to the manifolds $\bar{X}_1$ and $\bar{X}_2,$ and prove Theorem
\ref{inverseprob}  we proceed exactly
as section 8 of \cite{SA}.  The only new ingredient is  Lemma \ref{unique}, which, as mentioned before, is a consequence of
Lemma 2.3 of \cite{WV}.

\section{Appendix}
We give a short proof of the meromorphic extension of parabolically homogeneous distributions on $\hn=\rr_t\x\rr_z^{2n}$
\[u_\la(t,z)=(t^2+|z|^4)^{-\la}, \quad \Re(\la)<0\]
to $\la\in\cc$. It is very similar to the usual homogeneous distribution cases (see \cite[Th. 3.2.4]{HO}).
\begin{lem}\label{homog}
The family of parabolically homogenous distributions $u_\la$ on $\hn$ extends mermorphically to $\la\in\cc$ 
with only poles at each $\la_k=\frac{n+1}{2}+\demi k$ with $k\in\nn$, the residue of which is 
a distribution of order $2k$ supported at $0$.
\end{lem}
\textsl{Proof}: we consider the action of $u_\la$ again $f\in C_0^\infty(\hn)$, it is clear that
$\cjg \chi u_\la, f\cjd$ is analytic in $\cc$ for any $\chi\in C_0^\infty(\hn)$ such that $\chi=0$ near $0$. 
For the part near $0$, we use the parabolic coordinates
\[(R,u,\omega)\in (0,\infty)\x Q\to (R^2u,R\omega) \in\hn\setminus \{0\}, \quad Q:=\{u^2+|\omega|^4=1\}.\] 
If $v\in\rr^n\to \theta(v)\in S^{2n-1}$ is a parameterisation of the $2n-1$ dimensional sphere minus 
a point, then 
\[(u,v)\in(-1,1)\x \rr^n\to (u,\pm(1-u^2)^{\frac{1}{4}}\theta(v))\]
parameterizes each half of $Q$ and 
\[\psi_\pm:(R,u,v)\to (R^2u,\pm R(1-u^2)^\frac{1}{4}\theta(v))\]
can be used as changes of variable to compute 
\[\int_{\hn}(1-\chi(t,z))u_\la(t,z)f(t,z)dtdz.\] 
The Lebesgue measure pulls back to $R^{2n+1}G(u^2)dRdud\theta_{S^{2n+1}}(v)$ with 
$G(u):=(u+(1-u)^2)(1-u)^{\ndemi-1}$. 
Then taking $\chi$ with support in $\{R<1\}$ and depending only on $R$, we have to integrate
\[\int_{0}^1\int_{-1}^1\int_{S^{2n-1}} R^{-4\la+2n+1}(1-\chi(R))f\Big(R^2u, \pm R(1-u^2)^\frac{1}{4}\theta\Big)G(u^2)dRdud\theta_{S^{2n-1}}.\]
Then a Taylor expansion of $\psi_\pm^*f$ at $R=0$ gives for any $N\in\nn$ (using multi-index $\alpha$) 
\[f\Big(R^2u, R(1-u^2)^\frac{1}{4}\theta\Big)=\sum_{2i+|\alpha|<N} (\alpha!i!)^{-1}R^{2i+|\alpha|}u^i(1-u^2)^{\frac{|\alpha|}{4}}
\theta^\alpha\pl_t^i\pl^{\alpha}_zf(0,0)+O(R^{N})\]
and integration of $R^{-4\la+2n+1+2i+|\alpha|}(1-\chi(R))$ in $(0,1)$ extends meromophically to $\cc$
with pole at $\la=\frac{n+1+i+|\alpha|/2}{2}$, residue is clearly a distribution supported at $0$. 
Remark that for $|\alpha|$ odd, the residue involves the integral of $\theta^\alpha$ on $S^{2n+1}$,
which is easily seen to be $0$ by using change of variable $\theta\to-\theta$.
Note also that the residue of $\cjg u_\la,f\cjd$ at $\la_k$ is expressed
in terms of derivatives $\pl_z^{\alpha}f(0,0)$ with $|\alpha|=2k$ plus lower order 
derivatives, i.e. the principal term does not contain $\pl^{2k}_t$ derivatives.  
\qed\\

\end{document}